%Paper: alg-geom/9502013
%From: zjchen@fudan.ihep.ac.cn
%Date: Wed, 15 Feb 1995 19:34:49 BST-8
%Date (revised): Fri, 24 Mar 1995 19:14:02 BST-8

\input amstex
\magnification= \magstep1
\input amsppt.sty
\baselineskip15pt

\def\r{\rightarrow}
\def\o{\operatorname}
\def\ov{\over}
\def\e{\epsilon}
\topmatter
 \title Abelian automorphism groups of $3$-folds of general type
\endtitle
\author  Jin Xing Cai  \endauthor
\affil Department of Mathematics\\ East China Normal University \\
200062 Shanghai, China     \endaffil
\abstract
This paper is devoted to the study of abelian automorphism groups of
surfaces and $3$-folds of general type over complex number field $\Bbb C$.
We obtain a  linear bound in $K^3$ for
abelian automorphism groups of $3$-folds of general type
whose canonical divisor $K$ is numerically effective,
and we improve  on Xiao's results on   abelian  automorphism groups
of minimal smooth projective surfaces of general type.
More precisely,
the main results in this paper are the following.

{\bf Theorem 3.0.} Let $X$ be a  smooth 3-fold of general
type over the complex number field,
and  $K_X$ the canonical divisor of $X$.
Let $G$ be an abelian group of automorphisms of $X$ (i.e. $G\subset
\o{Aut} (X)$). Suppose $K$ is nef. Then
there exists a universal constant coefficient $c$ such that
 $$\# G \le c K_X^3 .$$
{\bf Theorems 6.5  and 7.1.}  Let $S$ be a minimal smooth surface of general
type over the complex number field,
and $K$ the canonical divisor of $S$.
Let $G$ be an abelian group of automorphisms of $S$ (i.e. $G\subset
\text{Aut} (S)$).
Suppose that $\chi(\Cal O_S)\ge 8$.
Then  $$\# G \le 36K^2 +24.$$
 Moreover, suppose that
 $K^2 \ge 181$,
and   that either $S$ has no pencil of curves of genus $g$, $3\le
g\le  5$ or $K_S^2\ge {12(g-1)\ov g+5}\chi(\Cal O_S)$
when  $S$ has a pencil of curves of genus $g$,
$3\le g\le  5$.
Then $$\# G \le 24K^2 +256.$$
\endabstract
\endtopmatter

\heading{ Introduction} \endheading

 Let $V$ be a  smooth projective variety  of general
type over the complex number field,
 $K_V$ the canonical divisor of $V$.
It is well-known that  the automorphism group $\o{Aut}(V) $ of $V$ is finite.
When $d:=\o{dim} V=1$,
a classical theorem of Hurwitz says that
$ \#\o{Aut(V)} \le 42\deg K_V$.
When $d=2$,
the automorphism groups  for a complex   smooth projective surface
  of general type
 have been thoroughly  studied by many authors.
In 1950,
Andreotti showed that the automorphism group $\text{Aut}(S)$ of
 a complex   smooth projective surface
$S$  of general type is finite and bounded by an
exponential function of $K^2$ \cite{A}.
Then Corti, Huckleberry and Saner independently proved that
$\vert \text{Aut}(S) \vert$  is bounded by a
polynomial  function of small degree in $K^2$  \cite{C}, \cite{HS}.
Recently,
progress has been achieved by Xiao,
who has obtained a linear bound for
$\text{Aut}(S)$,
in good analogy with the case of curves \cite{X2}, \cite{X3}:

{\it Xiao's theorem 1.
 Let $S$ be a minimal smooth surface of general
type over the complex number field,
 $K$ the canonical divisor of $S$.
Then $$\vert \text{Aut}(S) \vert  \le 42^2K^2 ,$$
and this bound is the best.}

Soon thereafter Chen has completely studied
the automorphism groups  for a  surface
  of general type with a pencil of genus $2$  \cite{Ch1}:

{\it Chen's theorem 1.
 Let $S$ be a minimal smooth surface of general
type with a minimal genus $2$ fibration,
and  $K$  as above.
Suppose  $K^2_S \ge 140$.
Then
(i) $\vert \text{Aut}(S) \vert  \le   504K^2 $;
moreover if f is not locally trivial,
then $\vert \text{Aut}(S) \vert  \le   288K^2 $.
And these are all the best bounds.}

When $d\ge 3$,
Xiao conjectured that under the assumption that $K_V$ is nef
(i.e., numerically effective),
then there is a number $\Delta_d$,
depending only on $d$,
such that $\#\o{Aut(V)} \le \Delta_d K_V^d$.

It is an intriguing problem to generalize these bounds to higher dimensions.
In the attempts to proving such a conjecture,
the order of abelian subgroups has a special importance
(see \cite{X1} for a discuss).

Let $G$ be an abelian subgroup of $\o{Aut}(V)$.
One has that $\#G  \le 2 \o{deg}K_V$\ $+8$ when $d=1$
\cite{N, W}.
For two dimensional  case,
Howard and Sommese proved that
 $\vert G\vert$ is bounded by the square of $K^2$
times a constant (among other things).
Soon before finishing the great work mentioned above,
Xiao considered the induced linear representation of
$G$ on the space of global sections of pluricanonical sheaves
and obtained a linear bound for $\#G  $
in terms of $K_V^2$ \cite{X1}:

{\it Xiao's theorem 2.
Let $S$, $K$ and $G$ be as above.
Suppose  $K^2_S\ge 140$. Then
 $\vert G \vert \le 52K_S^2+32$.}

For  abelian automorphism groups of   surfaces
  of general type with a pencil of genus $2$, Chen obtained
a best bound  \cite{Ch1}:

{\it  Chen's theorem 2.
 Let $S$ be a minimal smooth surface of general
type with a minimal genus $2$ fibration, and $G$ be as above.
Then  $\vert G \vert \le 12.5K_S^2+100$, if $K^2\ge 9$.}

The purpose of this paper
is to study the  abelian automorphism groups of   surfaces
and $3$-folds  of general type.
Of the two topics which we will concentrate on,
one is to
 give a $3$-dimensional
generalization of Xiao's results in \cite{X1} and
the other is to improve on Xiao's Theorem 2.
Our main results are the following.

{\it Theorem 3.0. Let $X$ be a  smooth 3-fold of general
type over the complex number field,
and  $K_X$ the canonical divisor of $X$.
Let $G$ be an abelian group of automorphisms of $X$ (i.e. $G\subset
\o{Aut} (X)$).
Suppose $K$ is nef.
Then
there exists a universal constant coefficient $c$ such that
 $$\# G \le c K_X^3 .$$}

{\it Theorems 6.5  and 7.1.  Let $S$ be a minimal smooth surface of general
type over the complex number field,
and $K$ the canonical divisor of $S$.
Let $G$ be an abelian group of automorphisms of $S$ (i.e. $G\subset
\text{Aut} (S)$).
Suppose that $\chi(\Cal O_S)\ge 8$.
Then  $$\# G \le 36K^2 +24.$$
 Moreover, suppose that
 $K^2 \ge 181$,
and   that either $S$ has no pencil of curves of genus $g$, $3\le
g\le  5$ or $K_S^2\ge {12(g-1)\ov g+5}\chi(\Cal O_S)$
when  $S$ has a pencil of curves of genus $g$,
$3\le g\le  5$.
Then $$\# G \le 24K^2 +256.$$}

The arguments here are inspired by the work  of Xiao  \cite{X1}.
We look at the induced linear representation
of $G$ on the space  $H_n$ of global sections of
pluricanonical sheaves $ \omega_V^{\otimes n} $.
Consider the natural map
$$ H_{n-t}  \otimes H_{n+t}  \oplus H_n  \otimes H_n  \rightarrow H_{2n}, $$
instead of$ H_n  \otimes H_n    \rightarrow H_{2n} $ as in \cite{X1},
we have that
under some weak conditions,
there are at least   two $G$- semi-invariants in  $H_{2n}$
corresponding to a same character of $G$.
Then  the corresponding
divisors  in $ \vert 2nK_S\vert$.
  generate a pencil $\Lambda$
whose general fiber $F$ is fixed by $G$.
Therefore $\#G$  is limited by the order of the group of automorphisms of the
normalization $\tilde F $ of $F$ as a smooth curve when
$\o{dim}(V)=1$ or that of
 the
minimal model $\tilde F_0 $ of the desingularization
  $\tilde F $ of $F$ as a minimal  smooth surface of general type
when $\o{dim}(V)=2$.
Theorems 3.0, 6.5 and 7.1 are obtained in this way.

I am convinced that essentially the same methods used to prove
Theorem 3.0 can be used to prove the statement for higher
dimensional algebraic varieties of general type.

I hope to return to this subject in a later occasion.

\heading Chapter I   Abelian automorphism groups of $3$-folds:
a linear bound \endheading

\heading \S 1 Preliminaries on $3$-folds \endheading

We work throughout over the complex number $\Bbb C$.

For the reader's convenience,
we collect some Preliminaries on $3$-folds of general type
in this section.
For Preliminaries on surfaces we refer the reader
to the standard text, e.g. \cite{BPV}.
We use the standard notation as in \cite{BPV}
unless otherwise stated.

Let $X$ be a nonsingular projective $3$-fold over $\Bbb C$,
and $D\in \o{Div} (X)$,
where $ \o{Div} (X)$ is a free abelian group generated by Weil
divisors on $X$.
We say that $D$ is nef if $D.C\ge 0$ for any curve $C$ on $X$,
and   that $D$ is big if the Iitaka dimension $\kappa (D, \ X)=3$ \cite{Ii}.
If $D$ is nef,
then for any integer $0\le i \le 3$ we know that $D^i.W\ge 0$
for every codimension $i$ subvariety $W$ of $X$ (cf.\cite{Ha, p.34}).
We say that $X$ is of general type if $K_X$ is big.
If $K_X$ is nef, then $K_X$ is big iff $K_X^3>0$ (Sommese,
cf.  \cite{Ka, Lemma 3}).

We denote the linear equivalence and the numerical equivalence
by $\equiv$ and $\sim $ respectively.
Denote by
$\Phi_n$ the rational map associated with  the complete linear system
$\vert nK_X \vert$.

By Hirzebruch-Riemann-Roch theorem we have
$$ \chi (\Cal O_X (nK_X))={1\over 12}(2n-1)n(n-1)K_X^3+(1-2n)\chi(\Cal O_X).
\tag 1.1$$

\proclaim{ Fact 1.2} (cf. \cite{Ka, Theorem 2 },  \cite{V, Theorem 1})
If $L$ is an nef and big line bundle on $X$,
then $H^i(X,\ L\otimes \omega_X)=0$ for $i>0$.
\endproclaim
\proclaim{ Corollary 1.3} If $X$ is of general type and $K_X$ is nef,
then

(i) $p_n:=h^0(nK_X)={1\over 12}(2n-1)n(n-1)K_X^3+(1-2n)\chi(\Cal O_X) $
for $n\ge 2$;

(ii) $\chi(\Cal O_X)\le K_X^3/6$, and $K_X^3$ is even;

(iii) $p_n\ge 5$ for $n\ge 3$.
\endproclaim
\proclaim{ Fact 1.4} Let $X$ be a smooth $3$-fold of general type and
 $f \colon X \r B$ be a fiber space,
where $B$ is a smooth projective curve.
By the {\it easy addition formula} (cf.  \cite{U, Chapter 6}),
the general fiber of $f$ is a smooth projective surface of general type.
\endproclaim
\proclaim{ Fact 1.5} (cf. \cite{Ma, Theorem 8})
Let $X$ be as in Corollary 1.3.
Then $\Phi_n$ is birational for $n\ge 8$.
\endproclaim
Finally we need a rough upper estimate of $-\chi (\Cal O_X)$
in terms of $K_X^3$.

\proclaim{ Proposition 1.6}
Let $X$ be as in Corollary 1.3.
Then $-\chi (\Cal O_X) \le {5 \over 2}K_X^3+1 $.
\endproclaim
\demo{Proof} We can suppose $-\chi (\Cal O_X)>0$;
consequently $p_2\ge 2$.
 Consider the following commutative diagram:
$$\alignat3 X' &   \overset f \to \longrightarrow  &  C & \hskip6cm\cr
\pi\downarrow & \searrow g &  \downarrow  & h \cr
X & \overset\Phi_2\to {\cdots \rightarrow}  &  W: &  =\Phi_2(X)
\subset \Bbb P^{p_2-1} \cr
 \endalignat$$
where
$\pi$ is  a succession of blowing-ups
with nonsingular centers such that $g:=\Phi_2 \circ \pi$ is a morphism,
and  $g=h\circ f$ is the stein factorization.
Let $\alpha=\o{deg}h$,
and $F$ be a general fiber of $f$.
Let $H$ be a hyperplane section of $W$  in $\Bbb P^{p_2-1}$, and
 $\beta$ the degree of $W$ in $\Bbb P^{p_2-1}$.
We have $$\pi^*(2K_X)\equiv f^*h^*(H)+Z,  $$
where $Z$ is the fixed part of $\vert \pi^*(2K_X)\vert$.

If $\o {dim} W=1$,
then $\pi^*(2K_X)\sim  \alpha\beta F+Z$.
Multiplying this equality by $\pi^*(K_X)^2$,
we have $$2K_X^3=  \alpha\beta \pi^*(K_X)^2. F+\pi^*(K_X)^2.Z.$$
Since $\pi^*(K_X)$  is nef and big,
and $F$ is nef and $F\not\sim 0$,
we have $\pi^*(K_X)^2.F\ge 1$ and $\pi^*(K_X)^2.Z\ge 0$.
Hence  $2K_X^3\ge \beta$.

If $\o {dim} W=2$,
then $\pi^*(2K_X)^2\sim  \beta F+\pi^*(2K_X).Z+(f^*h^*H).Z$.
Multiplying this equality by $\pi^*(K_X)$,
for the same reason as above
we have $$4K_X^3=  \beta \pi^*(K_X). F+2\pi^*(K_X)^2.Z
+\pi^*(K_X).f^*h^*H.Z\ge \beta.$$

If $\o {dim} W=3$,
for the same reason as above
we have
$$\pi^*(2K_X)^3-(f^*h^*H)^3=Z(\pi^*(2K_X)^2+\pi^*2K_X.f^*h^*H
+(f^*h^*H)^2)\ge 0.$$
Thus  $8K_X^3\ge \beta$.

Summing up the above inequalities,
we have  $8K_X^3\ge \beta$.

Now the result follows from the well-known fact that
$\beta\ge \o{codim} W+1$ (cf. e.g. \cite{Mu, p. 77}).
\qed \enddemo

\heading \S 2  Mid-points of finite lattice points:
some technical lemmas \endheading

Our goal in the next  section is to find a linear bound
for the order of abelian automorphism groups of $3$-folds of general type.
Our technique is to reduce the problem to that of  estimating the number of
 mid-points of some  finite lattice points in a linear space.

{\bf Definition 2.1} (\cite{X1}).
Let $\Cal  A$ and $\Cal  B$ be finite sets of points in a linear space $P$.
We define $\Cal  A.\Cal  B$  the set of mid-points ${1 \over2}(p+q)$
 of two points $p$
in $ \Cal  A$ and $q$ in $ \Cal  B$ ($p$ and $q$ may be the same point
if $\Cal  A \cap \Cal  B\not=\emptyset$; so   $\Cal  A \cap \Cal  B
\subset \Cal  A.\Cal  B$).

We define the {\it dimension} of $\Cal  A$ to be the
dimension of the  (affine) space generated by $\Cal  A$.
  Let $\Cal  A$ be a  finite set of points in  $P$,
 and $\Cal  B$ a subset of $\Cal A$.
The set  $\Cal  B$ is said to be {\it relatively convex} in  $\Cal  A$,
if no point of  $\Cal  A-  \Cal  B$ is contained in the convex hull
 of  $\Cal  B$.
The set  $\Cal  B$ is called {\it integrally convex} if it is relatively
convex in some lattice  $\Cal  A$ generating $P$.
With such a lattice  $\Cal  A$ fixed,
we will call the  points in  $\Cal  A$ integral points.

A {\it chain} in a set  $\Cal  B$ is by definition
a series of points $p_1$, ..., $p_n$ in $\Cal  B$
such that the vectors $p_i-p_{i-1}$ ($i=2$,...,$n$) are equal.
In this case,
$n$ is called the {\it length} of the chain.
If  $\Cal  B$ is integrally convex (in a fixed lattice  $\Cal  A$)
and $p$, $q$ are two points in  $\Cal  B$,
then the integral points on the line segment joining $p$ and $q$
form a chain in  $\Cal  B$ in an obvious way.

\proclaim{Lemma 2.2}  Let $ \Cal A_i$
( $ i=1,\ 2,\ 3$) be  finite integral sets in a $\Bbb Q$-linear space
with $  \Cal A_1 \subset   \Cal A_2  \subset   \Cal A_3$.
Let $P_i$ be the enveloping space of $\Cal A_i$ for $i=1$, $2$, $3$.
Assume $\text{dim}\  P_2<\text{dim}\  P$.
Then there exists an integral linear map $ \varphi  $
(i.e., $\varphi$ maps integral points to integral points)
from $P_3 $ to $P_2$ such that

(i) $\varphi|_ {\Cal A_3}$  is injective, and $ \varphi|_{P_2}$
 is identity;

(ii) $ \#(  \Cal A_1.  \Cal A_3 \cup   \Cal A_2.  \Cal A_2)
 \ge  \#(  \Cal A_1.\varphi(  \Cal A_3) \cup   \Cal A_2.  \Cal A_2).
$  \endproclaim
\demo{Proof}
Let $x_1$,\dots, $x_n$ be a (integral) basis of $P$.
We can suppose $\text{dim}\  P_2=\text{dim}\  P-1$,
and $P_2$ is the hyperplane $x_n=0$.
Since $\Cal A_3 $ is a finite set,
we can choose an integer $t$
such that $\vert a_i\vert+\vert b_i\vert <t$
($i=1,\ ...,\ n$) for any two points $(a_1,\ ...,\ a_n)$, $(b_1,\ ...,\ b_n)$
in $\Cal A_3$.
We define
$$\align \varphi :\ \ \  P & \rightarrow P_2 \cr
(x_1,\ x_2,\ ...,\ x_n) &
\mapsto (x_1+tx_n,\ x_2+tx_n,\ ...,\  x_{n-1}+t x_n,\ 0). \endalign$$
Clearly, $\varphi$ satisfies (i) by the choice of $t$.
Now the mid-point of two points $p$ in $\Cal A_1$ and $q$ in
$\varphi(  \Cal A_3)$ is the image of
the mid-point of two points $p$ in $\Cal A_1$ and $\varphi^{-1}(q)$ in
$  \Cal A_3$,
and if ${1\over 2}(p_1+q_1)={1\over 2}(p_2+q_2) $
for  $p_i$ in $\Cal A_1$ and $q_i$ in
$  \Cal A_3$, then
${1\over 2}(p_1+\varphi(q_1))={1\over 2}(p_2+\varphi(q_2))$.
Hence  $\varphi$ satisfies (ii).\qed \enddemo

Similarly, we have the following two lemmas.

\proclaim{Lemma 2.2.1}  Let $ \Cal A_i$
,  $P_i$  be as in Lemma 2.2.
Assume $\text{dim}\  P_1<\text{dim}\  P_2$, and $P_2=P_3$.
Then there exists an integral linear map $ \varphi  $
(i.e., $\varphi$ maps integral points to integral points)
from $P_2 $ to $P_1$ such that

(i) $\varphi|_ {\Cal A_3}$  is injective, and $ \varphi|_{P_1}$
 is identity;

(ii) $ \#(  \Cal A_1.  \Cal A_3 \cup   \Cal A_2.  \Cal A_2)
 \ge  \#(  \Cal A_1.\varphi(  \Cal A_3) \cup  \varphi( \Cal A_2).
 \varphi( \Cal A_2)).
$  \endproclaim

\proclaim{Lemma 2.2.2}  Let $ \Cal A_i$
,  $P_i$ be  as in Lemma 2.2.
Assume  $P:=P_1=P_2=P_3$ , and the  dimension of $P$ is $d\ge 2$.
Let $x_1$,\dots, $x_n$ be a (integral) basis of $P$.
Then there exists an integral linear map $ \varphi  $
(i.e., $\varphi$ maps integral points to integral points)
from $P $ to  the hyperplane  $P'\colon x_n=0$,  such that
{}.
Then there exists an integral linear map $ \varphi  $
(i.e., $\varphi$ maps integral points to integral points)
from $P $ to $P'$ such that

(i) $\varphi|_ {\Cal A_3}$  is injective; and

(ii) $ \#(  \Cal A_1.  \Cal A_3 \cup   \Cal A_2.  \Cal A_2)
 \ge  \#( \varphi( \Cal A_1).\varphi(  \Cal A_3) \cup  \varphi( \Cal A_2).
 \varphi( \Cal A_2)).
$  \endproclaim

{\bf  Definition 2.3} ( \cite{X1}).
Let $\Cal  A$ be a finite integral set in a space $P$
with a basis $x_1$, ..., $x_n$.
  The {\it arrangement}  of
$ \Cal  A $ with respect to $x_i$ is a new integral set $\Cal  A'$
such that a point $(a_1,\ .\ .\ .,\ a_n)$ is contained in $\Cal  A'$
iff

(i) $a_i\ge 0$, and

(ii) there are at least $a_i+1$ points  $(b_1,\ .\ .\ .,\ b_n)$
in $\Cal  A$ such that $a_j=b_j$ for all $j\not=i$.

It is immediate that $\Cal  A'$ is integral,
and has the same cardinality and dimension as that of  $\Cal  A$.
But if $\Cal  A$ is convex,
 $\Cal  A'$ need not be convex.
Clearly, If  $\Cal  A_i$ ($i=1,\ 2,\ 3$) are  finite integral sets in  $P$
with  $\Cal  A_1 \subset \Cal  A_2 \subset \Cal  A_3$,
then by the definition we
have   $\Cal  A'_1 \subset \Cal  A'_2 \subset \Cal  A'_3$.

 \proclaim{Lemma 2.4} Let $ \Cal A_i$ ($i=1$, $2$, $3$)
be  finite integral sets of dimension $d\ge 3$
with $  \Cal A_1 \subset   \Cal A_2  \subset   \Cal A_3$,
and let $ \Cal A'_i$ be  the arrangement of
$ \Cal A_i$ with respect to a coordinate  axis.
Then
$$ \#(  \Cal A_1.  \Cal A_3 \cup   \Cal A_2.  \Cal A_2 )
 \ge  \#(  \Cal A'_1.  \Cal A'_3 \cup   \Cal A'_2.  \Cal A'_2).
$$
\endproclaim
 \demo{Proof} We prove the lemma for the case $d=3$;
for the case $d>3$, the proof is similar, and left to the reader.

  For the sake of simplicity,
let $i=1$,
and use
 $x$, $y$, $z$ instead of $x_1$, $x_2$, $x_3$.
Take any two integers $y_m$, $z_m$,
and assume that there are $k$ points in
$ \Cal A_1'. \Cal A_3' \cup  \Cal A_2'. \Cal A_2'$ with $y=y_m$,
$z=z_m$.
Then there exist two points $p=(x_1, \ y_1 ,\ z_1)$ and
$q=(x_2, \ y_2 ,\ z_2)$ in
$ \Cal A_3'$ whose mid-point is $ ({1 \over2}(k-1), \ y_m, \ z_m)$.
Now  $ ({1 \over2}(k-1), \ y_m, \ z_m)$ is either in
$ \Cal A_1'. \Cal A_3'$ or in $ \Cal A_2'. \Cal A_2'$.
Hence we have either $p  \in  \Cal A_1'$,  $q  \in  \Cal A_3'$ or
 $p, \ q  \in  \Cal A_2'$.
Now we
suppose  $p  \in  \Cal A_1'$,  $q  \in  \Cal A_3'$. (for the latter case,
the proof is similar.)
By the definition of arrangement, this means that there
are at least $x_1+1$ points with
$y=y_1$, $z=z_1$ in $\Cal A_1$
and  at least $x_2+1$ points with
$y=y_2$, $z=z_2$ in $\Cal A_3$.
Because $x_1+x_2+1=k$,
we see that the points in  $\Cal A_1$ with
$y=y_1$, $z=z_1$ and the  points  in $\Cal A_3$  with
$y=y_2$, $z=z_2$  produce at least $k$
mid-points   with
$y=y_m$, $z=z_m$ in $\Cal A_1. \Cal A_3$:
in fact,
let the $x_1+1$ (resp. $x_2+1$) points in the first (resp. second)
row of  $\Cal A_1$ (resp.  $\Cal A_3$)
be $p_1$, ..., $p_{x_1+1}$
(resp.  $q_1$, ..., $q_{x_2+1}$) such that if $i<j$
then the $x$ coordinate of $p_i$ (resp. $q_i$)
is less than  $p_j$ (resp. $q_j$).
Then the mid-points
$$ {1\over 2}(p_1+q_1), \  {1\over 2}(p_1+q_2), \ ..., \
 {1\over 2}(p_1+q_{x_2+1}), \
{1\over 2}(p_2+q_{x_2+1}),\ ...\ , \  {1\over 2}(p_{x_1+1}+q_{x_2+1})
 $$
form the desired subset of  $\Cal A_1 .\Cal A_3$ \qed  \enddemo

\proclaim{Lemma 2.5}  Let $ \Cal A_i$
( $ i=1,\ 2,\ 3$) be  finite integrally convex sets
with $  \Cal A_1 \subset   \Cal A_2  \subset   \Cal A_3$,
and the dimension of $\Cal  A_1 $ is $3$.
Suppose that the length of the longest chain
 in $\Cal  A_3$ (resp.  $\Cal  A_2$)
is less then ${1\over 6}\# \Cal  A_3$ (resp. ${1\over 4}\# \Cal  A_2$).
And suppose
$$\# \Cal  A_2\ge 21,\ \
\# \Cal A_3 \le 2 \# \Cal  A_2. $$
Then $$ \# ( \Cal  A_1. \Cal  A_3\cup  \Cal  A_2. \Cal  A_2)
\ge \min \cases\# \Cal  A_3+3\# \Cal  A_2-23; \cr
{5\over 6}\# \Cal  A_3+{10\over 3}\# \Cal  A_2-10;\cr
{5\over 6}\# \Cal  A_3+{13\over 4}\# \Cal  A_2-2;\cr
{7\over 1 2}\# \Cal  A_3+{15\over 4}\# \Cal  A_2-6;.\cr
{1\over 2}\# \Cal  A_3+4\# \Cal  A_2-4; \cr
5\# \Cal  A_2-31. \endcases$$
\endproclaim
 \demo{ Proof }  By  Lemmas 2.2, 2.2.1, 2.2.2,
we can suppose the dimension of  the enveloping space of $\Cal  A_3$ is $3$.
Let $p_1$, ..., $p_l$ be a longest chain in $ \Cal  A_2$.
Assume $p_1$ is the origin of the enveloping space, and that
$p_2=(1, \ 0, \ 0)$.
Let $x$, $y$, $z$ be this basis. Then
$$   \root 3\of {\#\Cal  A_2}\le l
 \le {1 \over 4}  \#  \Cal  A_2.$$
  Arrange
$ \Cal  A_n$ ($i=1,\   ...,\   3$) with respect to $x$-axis,
then with respect to $y$,
and then to $z$.
By Lemma 4 ,
we have only to count the number of points in the set of mid-points of
the new set $ \#(  \Cal A'_1.  \Cal A'_3 \cup   \Cal A'_2.  \Cal A'_2) $
thus produced.

Let $m_{xi}$ (resp. $m_{yi}$,  $m_{zi}$ )
be the number of points of $ \Cal  A'_i$ on the $x$ (resp. $y$, $z$)
axis,
$i=1,\ 2, \ 3$.
We have
$m_{x2}=l \le { 1 \over 4} \# \Cal  A_2$
(cf. \cite{5, p. 625}) and $m_{x3} \le {1 \over 6} \# \Cal  A_3$,
$ m_{x1} \le m_{x2} \le m_{x3}$, etc.
If the basis elements $y$ and $z$ are chosen carefully,
$ \Cal  A'_1$ will also be of dimension $3$.
We may assume that the points of $ \Cal  A'_3$
have either $z=0$, or $z=1, \ y=0$
 (cf.  \cite{5, p. 626}).

Now for  simplicity of the notation,
we replace $ \Cal  A_i'$ by $ \Cal  A_i$.
Remark that   $(0,\ 0,\ 1)$ is in $\Cal  A_1$,
and although $ \Cal  A_i$ is not convex now,
it has the property that if a point $(a, \ b, \ c )$ is in $ \Cal  A_i$,
then all integral points  $(a', \ b', \ c' )$ such that
$$0 \le a' \le a,\ \  0 \le b' \le b,\ \  0 \le c' \le c$$
are in $ \Cal  A_i$ too.

Let $ \Cal  A_{i0}$ be just the subset of points of $ \Cal  A_i$
in the plane $z=0$.
Let $t_i$ ($i=1,\ 2, \ 3$) be the number of points in $ \Cal  A_i$ with
$z=1$.
Clearly,
$$1 \le t_1 \le t_2 \le t_3 \le m_{x3} \le {1 \over 6}
 \# \Cal  A_3.$$

In what follows we denote by $\Cal  S$ the set
$ \Cal A_1.  \Cal A_3 \cup   \Cal A_2.  \Cal A_2$.

Consider the mid-points  ${1 \over2}(r+s)$
with $r,\ s  \in  \Cal  A_ 2$ and $r$ on the $x$-axis and $s$
on the $y$-axis.
Two of  such points     ${1 \over2}(r_1+s_1)$ and  ${1 \over2}(r_2+s_2)$
are different if $r_1\not= r_2$ or $s_1\not= s_2$.
Consider the mid-points  ${1 \over2}(p+q)$
with $p=(0, \ 0, \ 1) \in  \Cal  A_ 1$ and $q$  is in the plane  $z=0$.
Therefore  $ \Cal S$ contains at least $m_{x2} m_{y2} + \#
\Cal  A_{30}$ points.
Hence we can assume $m_{x2} m_{y2}  \le {13\over 4} \# \Cal  A_2$.
Consequently,
$$ m_{x2}+m_{y2} \le {1 \over 4}  \#  \Cal  A_2 +12, \tag 2.5.1$$
as we have
$  \root 3\of {\#\Cal  A_2}\le m_{x2}
 \le {1 \over 4}  \#  \Cal  A_2$,   $m_{y2}\le {1 \over 4}  \#  \Cal  A_2$.

First, we have
$$\align  \#  ( \Cal A_{20}.\Cal A_{20}) &
=4\# \Cal  A_{20} -2(m_{x2}+m_{y2})+1, \tag 2.5.2\cr
 \#  \Cal S \ge \# ( \Cal A_2.\Cal A_2) &
=(t_2-1)(m_{y2}-3)+5\# \Cal  A_2 -2(m_{x2}+m_{y2})-3 .
\tag 2.5.3 \endalign $$
See  \cite{X1,  Lemma 6} for proofs of $(2.5.2)$ and $(2.5.3)$.

Second,
we estimate $\# \Cal  S$ by
considering not only the set
$ \Cal A_2.  \Cal A_2$ but also the set $   \Cal A_1.  \Cal A_3$.

Since  $p:= (0,\ 0,\ 1)$ is in $\Cal  A_1$,
 the mid-points  in
$$\Cal  B_1:=\bigg\{{1 \over2}(p+q) , q \  \text{is in } \Cal  A_{30} \bigg\}$$
are in $ \Cal A_1.  \Cal A_3$;
clearly, the  mid-points  in  $ \Cal B_2:= \Cal  A_{30}.\Cal  A_{30} $ and
$$\Cal  B_3:=\bigg\{{1 \over2}((a, \ 0,\ 1)+(b,\ 0,\ 1)) , 0\le a,\ b\le t_2-1,
a,\ b \in \Bbb Z \bigg\} $$
are in $ \Cal A_2.  \Cal A_2$.
It is easy to see that  $ \Cal B_i\cap  \Cal B_j=\emptyset$ for $i\not= j$.
Let
$ \Cal B= \Cal B_1\cup \Cal B_2\cup \Cal B_3. $
Then
$$  \# \Cal B=
 \#\Cal  A_{30}+\#( \Cal  A_{20}. \Cal  A_{20})+2t_2-1. $$
By ($2.5.1$) and ($2.5.2$), we have
$$\align  \# \Cal  S \ge & \# \Cal B
= \#\Cal  A_3+4\# \Cal  A_2-2(m_{x2}+m_{y2})-t_3-2t_2 \tag 2.5.4
\cr \ge & \# \Cal  A_3+{7\over 2}\# \Cal  A_2-t_3-2t_2-24.\tag 2.5.5
 \endalign $$

Now we consider the following cases separately.

 {\it Case\/} I1. ${1\over 2 } t_3 <t_2 \le {1 \over 8} \# \Cal  A_2$.
 \ \ By ($2.5.5$), we have
$$  \# \Cal  S \ge
 \# \Cal  A_3+{7\over 2}\# \Cal  A_2-4t_2-23
\ge  \# \Cal  A_3+3\# \Cal  A_2-23 . $$

   {\it Case\/} I2. $t_2\le \text{min}
\{{1\over 2 } t_3 ,\  {1 \over 8} \# \Cal  A_2$\}.
\ \  Let $a=t_3-2t_2+1$.
Since $p:=(0,\ 0,\ 1)$ is in $\Cal  A_1$,
 the mid-points  in
$$\Cal  C:=\bigg\{{1 \over2}(p+(t_3-1,\ 0,\ 1)) ,
{1 \over2}(p+(t_3-2,\ 0,\ 1)) ,\ ...,\ {1 \over2}(p+(t_3-a,\ 0,\ 1)) \bigg\}$$
are in $\Cal  S$;
it is clear that $ \Cal B\cap  \Cal C=\emptyset$,
and  $\# \Cal  C=a=t_3-2t_2+1$.
Now by ($2.5.5$), we have
$$  \# \Cal  S \ge
 \# \Cal  B+\# \Cal  C
 \ge  \# \Cal  A_3+3\# \Cal  A_2-23 . $$

   {\it Case\/} II1. $t_2 \ge {1 \over 8} \# \Cal  A_2 $
 and $ m_{y2} \ge 7$.
\ \  By ($2.5.1$) and ($2.5.3$), we have
$$  \# \Cal  S \ge
 4(t_2-1)+{9\over 2}\# \Cal  A_2-27
 \ge 5\# \Cal  A_2-31 . $$

  {\it Case\/} II2. $t_2 \ge {1 \over 8} \# \Cal  A_2 $
 and $ m_{y2} \le 6$.
\ \  Note  that we have $m_{y2}  \ge 3$
since  $ m_{x2} \le {1\over 4}\# \Cal  A_2$.
Let  $  e=  m_{y2} m_{x2} +t_2- \# \Cal  A_2 \ge 0,$
and let $n_i$ be the number of points
of $ \Cal  A_2 $ with $y=i, \ z=0$.
Clearly, $m_{x2}=n_0$.
 The mid-points in
$$ \align \Cal E:
=\bigg\{ &  {1 \over2} ((n_i-1, \ i, \ 0)+ (n_{i+1}-1, \ i+1, \ 0))
, \ ..., \cr & {1 \over2}( (n_{i+1}, \ i, \ 0)+ (n_{i+1}-1, \ i+1, \ 0)),\cr &
 \text{ provided} \   n_i>n_{i+1},\  i=0, \ ..., \ 3. \bigg\} \endalign$$
are in $\Cal S $,
 $ \Cal B \cap \Cal E=\emptyset$, and $\#  \Cal E=e$.
By ($2.5.4$),
we have
$$  \# \Cal  S \ge  \# \Cal B +\# \Cal  E
\ge  \#\Cal  A_3+4\# \Cal  A_2-2(m_{x2}+m_{y2})-t_3-2t_2+e.\tag 2.5.6
 $$

{ \it  Case\/} II2.1. $ m_{y2}=5$  or $ 6$.
\ \  If $ m_{y2}=6 $, then $ m_{x2}\ge{1\over 7}(\#\Cal  A_2+e)$.
By ($2.5.6$), we have
$$\align  \# \Cal  S \ge
 & \#\Cal  A_3+2\# \Cal  A_2+10m_{x2}-t_3-e-12
\cr \ge & {5\over 6}\# \Cal  A_3+{24\over 7}\# \Cal  A_2-12
 \ge {5\over 6}\# \Cal  A_3+{10\over 3}\# \Cal  A_2-10.
\endalign $$
Similarly, if  $ m_{y2}=5 $,
then we have
$  \# \Cal  S \ge
  {5\over 6}\# \Cal  A_3+{10\over 3}\# \Cal  A_2-10.$

{\it Case\/} II2.2.  $ m_{y2}=4$.
\ \  By ($2.5.6$), we have
$$ \# \Cal  S \ge
{5\over 6}\# \Cal  A_3+{13\over 4}\# \Cal  A_2-2+{e\over 4},$$
provided $m_{x2}  \ge{ 5 \over 24}
(\# \Cal  A_2+e)+1$.
This allows us to assume  $${1 \over 5} (\# \Cal  A_2+e) \le m_{x2}
 \le{ 5 \over 24}(\# \Cal  A_2+e).$$
Therefore, $t_2 \ge {1 \over 6}(\# \Cal  A_2+e)$.
Let $w$ be  the  integer part of
$ {1 \over 12}\# \Cal  A_3$.
Then the point $p=(w-1, \ 0, \ 1)$ is in $ \Cal  A_2$ since
$\# \Cal A_3 \le 2 \# \Cal A_2$ by the hypothesis,
and the mid-points in
$$ \Cal F: =\bigg\{{1 \over2}(p+(a, \ b, \ 0)),\
a\in \Bbb Z,\ a \ge w-1, \ b=0, \ ..., \ 3,
(a, \ b, \ 0) \ \text{is in}\   \Cal  A_{20}\bigg\}$$
are  in $\Cal S $,  $  \Cal B \cap \Cal F=\emptyset$,
 and $$\#  \Cal F=4(m_{x2}-w+1)-e
\ge 4(m_{x2}- {1 \over 12}\# \Cal  A_3+1)
-e.$$
Now by ($2.5.6$),
we have
$$\align  \# \Cal  S
\ge & \#\Cal  A_3+4\# \Cal  A_2-2(m_{x2}+m_{y2})-t_3-2t_2+e+ \# \Cal  F.
\cr \ge &  {1\over 2}\# \Cal  A_3+2\# \Cal  A_2+10 m_{x2}-2e-4
 \ge  {1\over 2}\# \Cal  A_3+ 4\# \Cal  A_2- 4.\endalign $$

{ \it Case\/} II2.3 $ m_{y2}=3$.
\ \  In this case  $   m_{x2}= t_2={1 \over 4} \# \Cal  A_2$.
The mid-points in
$$ \align \Cal G: = & \bigg\{{1 \over2}((a, \ 0, \ 1)+(a, \ b, \ 0)) ,\  \cr
& \text{where}\  { 1 \over12}\Cal \#  A_3  \le a \le
{1 \over4}
\# \Cal  A_2,\ a\in \Bbb Z, \  b=0, \ 1, \ 2.\bigg\}\endalign$$
are in $\Cal S$,
 $ \Cal B \cap \Cal G=\emptyset$ since
$2a>m_{x2}$, and $\#  \Cal G=3({1 \over 4}\# \Cal  A_2-{1 \over 12}\#
\Cal  A_3)$.
Then  by ($2.5.5$),
we have
$$  \# \Cal  S \ge  \#\Cal  B+\# \Cal  G
\ge { 7\over 12} \# \Cal  A_3+{15\over 4}\# \Cal  A_2-6. $$
Summing up above inequalities,
we get what we wanted.
  \qed  \enddemo

\proclaim{Lemma 2.6}   Let $ \Cal A_i$
( $ i=1,\ 2,\ 3$) be  finite integrally convex sets
with $  \Cal A_1 \subset   \Cal A_2  \subset   \Cal A_3$,
and the dimension of $\Cal  A_1 $ is  at least $4$.
Set $\e={1\ov 530}$.
Suppose that the length of the longest chain
 in $\Cal  A_3$
is less then $\e \# \Cal  A_3$, and $4 \# ( \Cal  A_2\ge \# \Cal  A_3$.
Then $$ \# ( \Cal  A_1. \Cal  A_3\cup  \Cal  A_2. \Cal  A_2)
\ge (1-\e)\# \Cal  A_3+{14(1-4\e)\ov 3}\# \Cal  A_2-57.$$
\endproclaim
 \demo{ Proof }  By  Lemmas 2.2, 2.2.1, 2.2.2,
we can suppose the dimension of  the enveloping space of $\Cal  A_i$ is $4$
for $i=1$, $2$, $3$.
Let $p_1$, ..., $p_l$ be a longest chain in $ \Cal  A_1$.
Assume $p_1$ is the origin of the enveloping space, and that
$p_2=(1, \ 0, \ 0, \ 0)$.
Let $x$, $y$, $z$, $u$ be this basis.
  Arrange
$ \Cal  A_i$ ($i=1,\   ...,\   3$) with respect to $x$-axis,
then with respect to $y$, to $z$,
and then to $u$.
By lemma 2.4,
we have only to count the number of points in the set of mid-points of
the new set $ \#(  \Cal A'_1.  \Cal A'_3 \cup   \Cal A'_2.  \Cal A'_2) $
thus produced.

Let $m_{xi}$
be the number of points of $ \Cal  A'_i$ on the $x$-axis,
$i=1,\ 2, \ 3$.
We have
$m_{xi} \le \e \# \Cal  A_3$.

If the basis elements $y$, $z$ and $u$ are chosen carefully,
then $ \Cal  A'_1$ will also be of dimension $4$.
We may assume that the points of $ \Cal  A'_3$
have either $u=0$, or $u=1, \ y=0,\ z=0$ as in the proof of
 \cite{X1, Lemma 6}.

Now for the simplicity of the notations,
we replace $ \Cal  A_i'$ by $ \Cal  A_i$.
Remark that although $ \Cal  A_i$ is not convex now,
it has the property that if a point $(a, \ b, \ c,\ d )$ is in $ \Cal  A_i$,
then all integral points  $(a', \ b', \ c',\ d' )$ such that
$$0 \le a' \le a,\ \  0 \le b' \le b,\ \  0 \le c' \le c, \  0 \le d' \le d$$
are in $ \Cal  A_i$ too.

Let
$ \Cal  A_{i0}$ be  the subset of points of $ \Cal  A_i$
in the plane $u=0$.
Let $t_i$ ($i=1,\ 2, \ 3$) be the number of points in $ \Cal  A_i$ with
$u=1$.
Clearly,
$$1 \le t_1 \le t_2 \le t_3 \le m_{x3} \le \e \# \Cal  A_3.$$

In what follows we denote by $\Cal  S$ the set
$ \Cal A_1.  \Cal A_3 \cup   \Cal A_2.  \Cal A_2$.

Consider the subset  $ \Cal  A_ {20}$.
Clearly it is of dimension $3$.
Since $t_2\le \e\#\Cal A_3\le 4\e \# \Cal A_2$,
we have $\# \Cal A_{20}=\# \Cal A_2-t_2\ge (1-4\e)\#  \Cal A_2$,
i.e., $\# \Cal A_2\le {1\ov 1-4\e}\# \Cal A_{20}$.

Let $k$ be the length of the longest chain in  $\Cal A_{20}$,
then $$k\le \e \#  \Cal A_3\le 4\e\#  \Cal A_2\le
 {4\e\ov 1-4\e}\#  \Cal A_{20}< {1\ov 6}\#  \Cal A_{20}.$$
Hence the condition of \cite{X1, Lemma 6}
is satisfied for   $\Cal A_{20}$;
one has that
$$ \#  ( \Cal A_{20}.\Cal A_{20})
\ge {14\ov 3}\# \Cal  A_{20} -57. \tag 2.6.1 $$

Second,
we estimate $\# \Cal  S$ by
considering not only the set
$ \Cal A_2.  \Cal A_2$ but also the set $   \Cal A_1.  \Cal A_3$.

Since  $r:= (0,\ 0,\ 0,\ 1)$ is in $\Cal  A_1$,
 the mid-points  in
$$\Cal  B_1:=\{{1 \over2}(r+q) , q \  \text{is in } \Cal  A_{30} \}$$
are in $ \Cal A_1.  \Cal A_3$;
clearly, the  mid-points  in  $ \Cal B_2:= \Cal  A_{20}.\Cal  A_{20} $
are in $ \Cal A_2.  \Cal A_2$.
It is easy to see that  $ \Cal B_1\cap  \Cal B_2=\emptyset$.
Let
$ \Cal B= \Cal B_1\cup \Cal B_2. $
Then
by ($2.6.1$), we have
$$\align  \# \Cal  S \ge & \# \Cal B\cr
\ge & \#\Cal  A_3+{14\ov 3}\# \Cal  A_{20}-t_3-57 \cr
\ge &  (1-\e)\# \Cal  A_3+{14(1-4\e)\over 3}\# \Cal  A_2-57. \qed \endalign $$
  \enddemo

Similarly, we have

\proclaim{Lemma 2.7}   Let $ \Cal A_i$
( $ i=1,\ 2,\ 3$) be  finite integrally convex sets
with $  \Cal A_1 \subset   \Cal A_2  \subset   \Cal A_3$,
and the dimension of $\Cal  A_1 $ (resp. $\Cal  A_2 $ ) is $\ge 2$
(resp. $\ge 3$).
Suppose that the length of the longest chain
 in $\Cal  A_3$ (resp.  $\Cal  A_2$)
is less then ${1\over 10}\# \Cal  A_3$ (resp. ${1\over 5}\# \Cal  A_2$).
Then $$ \# ( \Cal  A_1. \Cal  A_3\cup  \Cal  A_2. \Cal  A_2)
\ge
{9\over 10}\# \Cal  A_3+{16\over 5}\# \Cal  A_2-30. $$
\endproclaim

\heading{\S 3.  A linear bound of
abelian automorphism groups of $3$-folds of general type }\endheading

In this section we prove the following

 \proclaim{Theorem 3.0} Let $X$ be a  smooth 3-fold of general
type over the complex number field,
and  $K_X$ the canonical divisor of $X$.
Let $G$ be an abelian group of automorphisms of $X$ (i.e. $G\subset
\o{Aut} (X)$).
Suppose $K$ is nef.
Then
there exists a universal constant coefficient $c$ such that
 $$\# G \le c K_X^3 .$$
  \endproclaim

{\it Remark.}
It is easy to construct examples with $\# G$ increasing linearly with
$K^3$: let $S$ be a surface of general type with an abelian  automorphism
group $G_1$ of order $12.5K_S^2+100$ \cite{Ch},
and let $C$ be a curve of genus $2$  with an abelian  automorphism
group $G_2$ of order $12$ \cite{N}. Let $X=S\times C$.
Let $p_1$, $p_2$ be the projections of $X$ onto the two factors.
The group $G$ generated by $p_1^*G_1$ and $p_2^*G_2$ is
 an abelian  automorphism group of $X$ of order $25K_X^3+1200$.

We can give explicit  estimates in Theorem 3.0 but these are far from being
the best possible.
An interesting question is,
roughly speaking:
what is asymptotically the best upper bound for $G$?

The arguments here are inspired by Xiao's work
 on abelian automorphism groups of surfaces of general type \cite{X1}.
We consider the natural action of the abelian group $G$
on the space $H_n=H^0(X,\ nK_X)$,
for a fixed positive integer $n$.
Because $G$ is finite abelian,
such an action is diagonalisable,
in other words $H_n$ has a basis consisting of semi-invariant vectors.
Consider two such semi-invariants  $v_1,\ \ v_2$ in $H_n$,
with
$$\sigma (v_i)=\alpha _i(\sigma)v_i \ \ \text{for}\ \ \sigma \in G,$$
where
$ \alpha_i$ are the corresponding characters of $G$.
Suppose  that the  two  semi-invariants  $v_1, v_2$
corresponding to a same character of $G$ (i.e., $\alpha_1=\alpha_2$),
 and let
 $D_1$ and $D_2$ be  the corresponding
divisors  in $ \vert nK_X\vert$.
Then
 $D_1$ and $D_2$ generate a pencil $\Lambda$
whose general fiber $F$ is fixed by $G$.
Therefore $\#G$  is limited by the order of the group of automorphisms of the
minimal model $\tilde F_0 $ of the desingularization
  $\tilde F $ of $F$ as a minimal  smooth surface of general type.
But $\#\operatorname{Aut}(\tilde F_0)$ increases proportionally with $K^3_X$,
as $p_g(\tilde F)$ so does.

We  consider the natural map
$$ H_n  \otimes H_{3n}  \oplus H_{2n}  \otimes H_{2n}  \rightarrow H_{4n}, $$
which is compatible with the above actions of $G$:
i.e.,  if   $v_1 \in H_n$, $v_2 \in H_{3n}$ (resp. $w_i \in H_{2n}$
$i=1,  \ 2$)
are  two  semi-invariants,
then $v_1 \otimes v_2$ (resp.  $w_1 \otimes w_2$)
is  semi-invariant in $H_{4n}$.
If there are more than $\o{dim}\  H_i$ semi-invariants in
$H_i$ for some $i\le {3n}$,
then  there are semi-invariants in
$H_i$ (therefore in $H_{4n}$) with the  same character,
and we are done.
So we may assume that  there are exactly  $\delta_i=\o{dim}\  H_i$
($i\le {3n}$) semi-invariants $v_j^i$ ($j=1,\ ...,\ \delta_i$) in
$H_i$,
corresponding to mutually different characters.
Each vector $v_j^i$
corresponds to a unique divisor $D_j^i$ in $\vert iK \vert$.
The relation $v_j^i \otimes v_k^l=c v_m^{r}\otimes v_n^{s}$
(where $c$ is a constant,
and $i+l=r+s=4n$) in $H_{4n}$ translates to a  relation
$$D_j^i+ D_k^l= D_m^{r}+D_n^{s} \tag *$$
between these divisors.

Fix a semi-invariant $u\in H_n$ when $H_n\not= 0$.
Then $u$ corresponds to a unique divisor $U$ in $\vert nK \vert$ which is fixed
by $G$.
 We can consider the finite set  $\Sigma _i$
  of points corresponding to $D_j^i$
in a certain divisorial space $P_i$ defined in  \cite{X1, $\S 1$},
and there are  natural embeddings:
$$ \Sigma _n \rightarrow  \Sigma _{2n} \rightarrow  \Sigma _{3n},  \ \ \ \
\text{and} \ \  \ \  P_{n}   \rightarrow
P_{2n}   \rightarrow P_{3n}$$
defined by $U$.
In such a setting,
a semi-invariant in $H_{4n}$ of the form  $v_j^i\otimes v_k^l$
corresponds naturally to the mid-point of two points in
$ \Sigma _{3n}$
corresponding to $D_j^i$ and $ D_k^l$,
and a relation of the form (*) means that the corresponding mid-points
coincide.

Denote by $ \Cal S_1$ the set of
mid-points of two points $p$, $q$ in $ \Sigma_ {3n}$
such that either $p$ is in $ \Sigma_{n}$ and $q$ is in $ \Sigma_{3n}$
 or $p$ and $q$
are in $ \Sigma_{2n}$.
Now the problem has been reduced to
that of  comparing the number of points in $ \Cal S_1$
and  the dimension of $H_{4n}$.
In this way we show that for $n$ large enough
 the number of points in $ \Cal S_1$
is larger than   the dimension of $H_{4n}$.

We fix a smooth complex projective $3$-fold of general type $X$
in the future,
and let $K$ be  the canonical divisor of $X$,
$H_n=H^0(X,\ nK)$, and
$\chi=\chi(\Cal O_X)$.
We also fix  an abelian group $G$ of automorphisms of $X$.

For the reader's convenience,
we recall some notation defined in \cite{X1}.

{\bf Definition} (\cite{X1}).
Let $v_1$, ..., $v_{\delta_n}$ be a basis of $H_n$ consisting  of
 semi-invariants for the natural action of $G$ on $H_n$,
and $D_1$, ..., $D_{\delta_n}$
 the divisors in $\vert n K \vert$ corresponding to these vectors,
where $\delta_n=\o{dim}\ H_n$.
We say that $H_n$ is  {\it uniquely decomposable}
(under the action of $G$) if the set $\{D_i\}$ is uniquely determined,
or equivalently if there are exactly $\delta_n$ different characters
 for the natural action of $G$ on $H_n$.

Fix the divisor $D_1$.
Denote by $P'_n$ the set
$$\align  \{ \Bbb Q \text{-divisors}\  D \text{ on}\  X \ \vert
&\text{ there is an} \   m \in   \Bbb Z^+ \cr &
\text{such that}\  mD \ \text{is linearly equivalent to}\  mD_1 \}.\endalign$$
Denote by $[D]$ the element of $P'_n$ corresponding to
the $\Bbb Q$-divisor $D$.
We define addition and scalar multiplication  as follows:
$$\align [ D] + [ D'] & =
[ D+ D'- D_1 ] ,\cr  c[ D ] & = [ cD+(1-c)D_1]
 ,\ \ c \in \Bbb Q .\endalign $$
Then $P_n'$ is a generally infinite dimensional linear space,
with $ [  D_1 ] $ as the origin.

The subset $I$ in $P_n'$ of points  corresponding to
integral divisors  linearly equivalent to $D_1$ is an additive subgroup,
and there is a set of generators of $I$ which form a basis of $P_n'$.
Under such a  basis,
$I$ is a subset of points with integral coordinates.

Denote by $P_n$  the finite dimensional subspace generated by the set
$$ \{[ D_1],\  [ D_2 ] ,\ .\ .\ .\ ,\
[ D_{\delta_n}] \}.$$
Let $\Sigma_n$ be the finite set
 in $P_n$ consisting of  points  corresponding to
effective  divisors in  $\vert n K \vert$ fixed by $G$.
Then $P_n$ ,
therefore $\Sigma_n$,
is uniquely determined up to  choices  of $D_1$ only if $H_n$ is
uniquely  decomposable.
We will call the set $\Sigma_n$ {\it  a basic set } in $P_n$.

Clearly, $P_n$ depends on the choice of $D_1$;
but if we replace $D_1$ by another divisor,
say $D_i$,
$P_n$ differs only by an integral translation.
Because $\#\Sigma_n$ and the number of middle points of
 $\Sigma_n$,  which are all the properties about
$\Sigma_n$  we use later,
 are  integral translation
invariants,
it dose not  matter which such $D_i$ is selected.
Also,
$\Sigma_n$ is determined up to   integral translations as above
iff $H_n$ has exactly $\delta_n$ semi-invariants,
and thus iff there are $\delta_n$ different characters
for the  action of $G$ on $H_n$.

Fix a semi-invariant $u \in H_n$ for the natural action of $G$ on $H_n$.
Let $U$ be the divisor in $\vert nK \vert$ corresponding to $u$.
We have natural maps:
$$\align & H_n \overset {\otimes u} \to  \rightarrow
H_{2n} \overset { \otimes u} \to  \rightarrow H_{3n},\cr
&  \vert  nK \vert  \overset {+U} \to  \rightarrow
 \vert 2nK \vert  \overset {+U} \to  \rightarrow  \vert  3nK \vert ,\cr
 \text{and} \ \ &   P'_n \overset {l'_1} \to  \rightarrow
P'_{2n} \overset {l'_2} \to  \rightarrow P'_{3n}.\endalign$$
If we take $[ nU ]$ to be the origin of $P'_n$,
then $l'_i (i=1, 2)$ are embeddings of linear spaces.
We will identify
 $P'_n$ and $P'_{2n}$ as subspaces of $P'_{3n}$ in this way in the future.
Since $U$ is fixed by $G$,
 $l'_i$ ($i=1,\ 2$) induce  natural embeddings:
$$ \Sigma _n \rightarrow  \Sigma _{2n} \rightarrow  \Sigma _{3n},  \ \ \ \
\text{and}\ \
 P_n \overset {l_1} \to  \rightarrow
P_{2n} \overset {l_2} \to  \rightarrow P_{3n}.$$
Let $v $ (resp. $w$) be a semi-invariant in $H_n$ (resp. $H_{3n}$), and
 $p$ (resp. $q$) the corresponding
point in $ \Sigma_n$ (resp. $ \Sigma_{3n}$).
Let $D$ be the divisor in $\vert 4nK\vert$
corresponding to the vector $v \otimes w$.
Then $l_2 ([ {1 \over2}D ])$  corresponds to a point in $P_{3n}$,
 which is just
the mid-point ${1 \over2}(l_2l_1(p)+q)$.

Similarly, let $v_1,\ v_2$ be two  semi-invariants in $H_{2n}$,
 and
 $p$, $q$ the corresponding
points in $ \Sigma_{2n}$.
Let $D$ be the divisor in $\vert 4nK\vert$
corresponding to the vector $v_1 \otimes v_2$.
Then $l_2 ([ {1 \over2}D ])$  corresponds to a point in $P_{3n}$,
 which is just
the mid-point ${1 \over2}(l_2(p)+l_2(q))$.

 \proclaim{ Lemma 3.1}  Let  $\Sigma_n $  be a basic set in $P_n$.
Then

 (i)  The set  $\Sigma _n$ is  integrally convex
with respect to the lattice  $\Cal  L$
consisting of points corresponding to divisors linearly equivalent to
$nK$, and
$\Sigma _{ni} $ is relatively convex in $\Sigma _{n(i+1)} $
for $i\ge 1$,
as we consider $\Sigma _{ni} $ as a subset of $\Sigma _{n(i+1)} $
in the above way.

(ii) Suppose $P_n$ is uniquely decomposable.
 If the dimension of $\Sigma_n$ is less then $3$,
then $\o{dim}\Phi_n(X)\le 2$.

(iii)  Suppose $P_n$ is uniquely decomposable.
If $\Phi_n$ is birational,
then the dimension of $\Sigma_n$ is at least $4$.
\endproclaim
\demo{ Proof}
 (i)  See \cite{X1, Lemma 3}.

(ii) Suppose the  dimension of $\Sigma_n$ is $2$.
By hypothesis,
we have $$\Sigma_n = \{[ D_1],\  [ D_2 ] ,\ .\ .\ .\ ,\
[ D_{\delta_n}] \}.$$
We can suppose $[D_1]$ is the origin  of $P_n$,
and we can find two points, say $[D_2]$ and $[D_3]$ in $\Sigma_n$
such that

(a) $[D_2]$ and $[D_3]$ generate the additive group $\Cal L$,
and

(b) for each $[D_j]$ ($j>3$), we have $[D_j]=a_j[D_2]+b_j[D_3]$
with $a_j,\ b_j \in \Bbb Z^+$.

 By (b) we have $$D_j+(a_j+b_j-1)D_1=a_jD_2+b_jD_3$$ for $j\ge 4$,
that is, $$ v_j \otimes v_1^{\otimes(a_j+b_j-1)}=\xi_j
 v_2^{a_j}\otimes v_3^{b_j},\ \ \ \text{for some }\ \  \xi_j \in \Bbb C,
\ \ j\ge 4. \tag 3.1.1$$
Let $\phi \colon X \cdots \r \Bbb P^2$ be the rational map defined by
$(v_1\ :\ v_2\ :\ v_3)$.
Now for general points $p$, $q\in X$
(here by general we mean that $v_1(p)\not =0$ and $v_1(q)\not =0$),
such that $\pi(p)=\pi (q)$ (i.e., $ v_i(p)=\lambda v_i(q)$ for some
$\lambda \in \Bbb C$, $i=1$, $2$, $3$),
by (3.1.1) we have $ v_j(p)=\lambda v_j(q)$ for
 $j=4$, $\cdots$, $\delta_n$.
Hence $\Phi_n$ ($=(v_1\ :\ v_2\ :\ v_3\ : \cdots \ : v_{\delta_n})$)
is factor through $\phi$.

(iii) Otherwise,
by (ii) we can suppose that the dimension of $\Sigma_n$ is  $3$.
Then  $\Phi_n$
is factor through $\phi   \colon X \cdots \r \Bbb P^3$,
where $\phi$ is similarly defined as in (ii). This is a contradiction.
 \qed \enddemo

\proclaim{Lemma 3.2}Let $\Sigma_n$ be a basic set in $H_n$.
Then the length of the longest chain in $\Sigma _n$ is less then
$   {12\over (2n-1)(n-2)}  \# \Sigma_n$.
 \endproclaim
\demo{Proof}
Let $k$ be  the length of the longest chain in $\Sigma _n$.
Then there is a pencil $\Lambda$
of surfaces on $X$ such that there is a divisor in
$\vert nK_X \vert$ containing $k$ times a fiber $F$ of  $\Lambda$ .
Hence $nK_X \equiv kF+D$ for some divisor $D\ge 0$;
we have $kK_X^2.F\le nK_X^3$ since $K_X$ is nef.

We claim $K_X^2F\ge 1$.
Indeed, let $\pi \colon X' \r X$ be a succession of blowing-ups
with nonsingular centers such that $\phi \circ\pi$ is a morphism,
where $\phi$ is a rational map induced by  $\Lambda$ .
We set $\pi^*F\equiv F'+E$,
where $F'$ consists of  some fibers of  $\phi \circ\pi$,
and $E$ is an exceptional divisor for $\pi$.
We have $$ K_X^2.F=\pi^*(K_X)^2\pi^*.F=\pi^*(K_X)^2.F'+\pi^*(K_X)^2.E$$
Since $\pi^*K_X$ is nef,
we have $\pi^*(K_X)^2.E\ge 0$,
and Since $\pi^*K_X$ is nef and big,
$F'$ is nef, and $F'\not\sim 0$,
we have $\pi^*(K_X)^2.F'\ge 1$.

Hence we have $k\le nK_X^3<{12\over (2n-1)(n-2)}p_n$
by Corollary 1.3.\qed \enddemo

 \proclaim{Proposition 3.3} Suppose $K_X$ is nef.
Then  there exists an integer $N$ such that
 $H_{4N}$ is not uniquely decomposable.
  \endproclaim
\demo{Proof}
By Lemma 3.1 and Lemma 3.2, we can choose $N_0\gg 0$
such that for $n\ge N_0$,
the condition of  Lemma 2.6 is satisfied for
$ \Sigma _n \subset  \Sigma _{2n} \subset  \Sigma _{3n},  $
and we can suppose $H_i$ is  uniquely decomposable for $i< 4n$,
for otherwise the corollary is trivially true.
By Corollary 1.3, we have

 $$\align  & \# ( \Sigma _n .  \Sigma _{3n} \cup
    \Sigma _{2n} .  \Sigma _{2n}) -\o{dim}\  H_{4n}\cr
& \ge  (1-\e)\# \Sigma_{3n}  +{14(1-4\e)\over 3}\# \Sigma_{2n}-57-
 \o{dim}\  H_{4n} \cr
 &= (1-\e)\text{dim}H_{3n}  +{14(1-4\e)\over 3} \text{dim}H_{2n}-57-
 \o{dim}\  H_{4n} \cr & =({1-529\e \ov 18}n^3+O(n^2))K_X^3
\cr & >0\ \ \  \text{for} \ n\ge N_0\ \ \text{and}\ \  n\gg 0.\endalign$$
Hence there exists an integer $N$ independent of the choice of $X$ such
that $ \# ( \Sigma _N .  \Sigma _{3N}\mathbreak \cup
    \Sigma _{2N} .  \Sigma _{2N}) >\o{dim}\  H_{4N}$,
i.e., there are more than $\o{dim}  H_{4N}$- semi-invariants in $H_{4N}$.
\qed \enddemo

\demo{Proof of theorem 3.0} Let $b=4N$. The Proposition 3.3
guarantees that there is a pencil
  $\Lambda$
in $\vert bK \vert$ each of  whose  elements is a fixed divisor by $G$.
We consider the following commutative diagram:
$$\alignat3 X' &   \overset f \to\longrightarrow  &  C &\hskip6cm \cr
\pi\downarrow & \searrow g &  \downarrow  & h \cr
X & \overset\phi\to {\cdots \rightarrow}  &  \Bbb P^1 &\cr
 \endalignat$$
where $\phi$ is the rational map induced by  $\Lambda$,
$\pi$ is  a succession of blowing-ups
with nonsingular centers such that $g:=\phi \circ \pi$ is a morphism,
and $g=h\circ f$ is the stein factorization.
Let $d=\o{deg}h$,
and $F$ a general fiber of $f$.
Clearly, $d\le bK_X^3$.
By Fact 1.4 $F$ is a surface of general type.
Let $F_0$ be the minimal model of $F$.

The automorphism group $G$ can be lifted to be an  automorphism group
 $G'$ of $X'$.
(Clearly $\# G'=\# G$;
for any element $\sigma$ in $G$,
denote by $\sigma'$ the corresponding element in $G'$.)
Let $H$ be the stabilizer of $F$,
then the index $H$ in $G'$ is at most $d$.
We have a natural group homomorphism
$\rho \colon H \r \o{Aut}(F_0)$.
We claim that $\o{Ker}\rho=1$.
Indeed, by the definition of $\pi$,
we can identify  $X'-\pi^{-1}( Y) $ with $ X-Y$,
where $Y$ is the base locus of the moving part of   $\Lambda$.
Let $\sigma' \in \o{Ker}\rho $.
 Since $F$ is general,
 $\sigma'  $ equals identity on $X'-\text{
Exceptional divisors of } \pi$.
Since $X'-\pi^{-1}(Y)=X-Y$,
we have
$\sigma$ equals identity on $X-Y$.
Hence $\sigma$ equals identity on $X$ since the codimension of $Y$ in
$X$ is $\ge 2$.

Since $H$ is an abelian automorphism group of $F_0$,
we have (cf. \cite{X1}, \cite{Ca}.)
$$\# H \le \cases  36K_{F_0}^2+24, \ \  &\text{ if} \ \
\chi(\Cal O_{F_0}) \ge 8; \cr
270K_{F_0}^2, &\text{ otherwise}.\cr\endcases $$
Consequently,
we have $$\# G\le d\cdot \# H \le
\cases 335p_g(F_0)d, \ \  & \text{if}\ \  p_g(F_0)\ge 34; \cr
 270\cdot 9 \cdot 34d \le  270\cdot 9 \cdot 34
\cdot bK_X^3, & \text{otherwise}.\cr\endcases$$

Since  $p_g(F_0)=p_g(F)$, Theorem A  follows from the following claim
and Proposition 1.6.

{\it Claim.}  $dp_g(F)\le p_{b+4}(X)$.

{\it Proof of the claim.}\ \
 We have $$\pi^*(bK_X)\equiv f^*h^*(p)+Z, \tag 2.6 $$
where $p$ is a closed point of $\Bbb P^1$,
and $Z$ is the fixed part of $\pi^*(bK_X)$.
For general $p\in \Bbb P^1$,
we have $D:=f^*h^*(p)=F_1+\cdots +F_d$ such that $F_i$ are general
fibers of $f$ and $F_i\not=F_j$ when $i\not= j$.

Since $p_3(X)>0$ by Corollary 1.3,
we have $$(K_{X'}+D)\vert _D<(K_{X'}+D+\pi^*3K_X )\vert _D, $$
so $$ dh^0(K_F)=h^0(K_{X'}+D\vert _D)\le h^0((K_{X'}+D+\pi^* 3K_X)\vert _D).$$
Consider the following exact sequence
$$\align 0\r  \Cal O_{X'}(K_{X'}+\pi^*(3K_X))
& \r \Cal O_{X'}(K_{X'}+D +\pi^*(3K_X))
\cr   \r & \Cal O_D(K_{X'}+D+\pi^*(3K_X))\r 0.\cr\endalign$$
Since $\pi^*(3K_X)$ is nef and big,
by Fact 1.2 we have $$h^1( \Cal O_{X'}(K_{X'}+\pi^*(3K_X)))=0,$$
so $$h^0( \Cal O_{X'}(K_{X'}+D+\pi^*(3K_X))\vert_D)\le
h^0( \Cal O_{X'}(K_{X'}+D+\pi^*(3K_X))).$$
since $\pi_* \Cal O_{X'}(K_{X'})= \Cal O_X(K_X)$ (cf. \cite{Mo, p. 280} ) and
$D<\pi^*(bK_X)$ (by (2.6)),
 we have $$ h^0(\Cal O_{X'}(K_{X'}+D+\pi^*(3K_X)))\le
h^0( \Cal O_{X'}(K_{X'}+\pi^*((b+3)K_X)))
= p_{b+4}(X).$$
Summing up the above inequalities, we get $dp_g(F)\le p_{b+4}(X)$.
\qed \enddemo

 \heading Chapter II   Abelian automorphism groups of surfaces:
an improvement on  Xiao's  results
 \endheading

\
\
  \heading \S 4.  Automorphisms  of non-hyperelliptic curves
\endheading

It is well-known that for a complex curve of genus $g\ge 2$,
its  abelian  automorphism group
is of order  $\le 4g+4$ \cite{N}.
In this section  we give
 a further  analysis to abelian automorphism groups for non-hyperelliptic
curves. Our main
results are  Theorems 4.4 and 4.10 which shall be used in the sequel.

We begin by establishing notations.

{\bf Notation 4.1}.
Let $C$ be a smooth projective curve of genus $g\ge 2$,
$G$ an abelian subgroup of $\text{Aut}(C)$ .
We have a finite abelian covering
$$\pi \colon  C \rightarrow X:=C/G .$$
Let $q_1$, ..., $q_k$ be the points over which $\pi$ is
ramified,
and $r_i$ the  ramification number of $\pi$ over $q_i$.
We assume that $r_1\ge r_2\ge ... \ge r_k$.
Choose a point $Q_i \in C$ lying above $q_i$, and put
$$G(q_i)=\{\sigma \in G\  \vert\  \sigma Q_i=Q_i\}.$$
Since $G$ is assumed to be abelian,
$G(q_i)$ is well defined;
we have $\vert G(q_i)\vert =r_i$.
Using Hurwitz formula to the morphism $\pi$,
we get
$${(2g-2)\ov\vert G \vert} =2g(X)-2+\sum_{i=1}^k (1-{1\ov r_i}).\tag4.2$$
For convenience, we denote ${(2g-2)\ov \vert G \vert}$ by $\tau$.

\proclaim{Lemma 4.3} Notations being as above.
Let $\vert G \vert =p_1^{t_1}.....p_l^{t_l}$,
where $p_i$s are prime numbers,
and let $m_j=\Pi _{1\le i\le k,\ i\not=j}r_i$ for $j=1$,
..., $k$.
Suppose $g(X)=0$.
Then we have

(i) $m_j$ is a multiple of $\vert G \vert$ for any $j$;

(ii) For each $i$,
there are at least two points $q_j$, ${q_j'}$ such that $p_i \vert r_j$,
$p_i \vert {r_j'}$;

(iii) $\vert G \vert$ is a multiple of the least common multiple of $r_i$'s.

(iv) If $G$ is cyclic,
then  $\vert G \vert$ equals the least common multiple of $r_i$'s,
and
 for each $i$,
there are at least two points $q_j$, ${q_j'}$ such that $p_i^{t_i} \vert r_j$,
$p_i^{t_i} \vert {r_j'}$;
\endproclaim

\demo{Proof} The Galois group $G$ of $\pi$ is an abelian quotient
$\pi _1(X-\{ q_1, ..., q_k\})$ which is generated by $\gamma_1$,
..., $\gamma_k$,
where $\gamma_i$ is a small loop around $p_i$.
Let $\bar \gamma_i$ be the image of $\gamma_i$ in $G$.
$r_i$ being the order of $\bar \gamma_i$.
Since $\bar \gamma_1...\bar \gamma_k=1$,
we  have $G$  is generated by,
say,
$\bar \gamma_2$,
..., $\bar \gamma_k$.
Hence (i) follows.
Since $r_i\vert \ \vert G \vert$ for each $i$,
we have (iii).
Now (ii) follows from (i) and (iii).

 If $G$ is cyclic,
then  the order of the generator of $G$ divides
 the least common multiple of $r_i$'s.
Hence by (iii) we have that
$\vert G \vert$ equals the least common multiple of $r_i$'s.
Since say, $r_1$ is the order of $\bar \gamma_1^{-1}=\bar \gamma_2
\dots \bar \gamma_k$,
we have $r_1$ divides
 the least common multiple of $r_i$'s.
Hence if
 $p_i^{t_i} \vert r_1$, then there is a $j\ge 2$ such that
$p_i^{t_i} \vert {r_j'}$. \qed \enddemo

{\bf Remark} The assumption of being cyclic is indispensable in
Lemma 4.3 (iv).
For example, let  $C$ is the plane Fermat curve $x^5+y^5+z^5=0$,
and $G$ consists of the automorphisms $(x, \ y, \ z)\rightarrow (\mu x,\ \nu
y,\ z)$,
where $\mu$, $\nu$ are $\text{exp}({2\pi i\ov 5})$.
Then $X/G\simeq\Bbb P^1$,
 $\vert G \vert =25$,
and $r_i=5$ for $i=1$, $2$, $3$.
\proclaim{Theorem 4.4} Notations being as in (4.1).
Suppose that $C$ is nonhyperelliptic.
Then $$\vert G \vert \le 3g+6$$
with one exception below.
The exceptional $C$ is the plane Fermat curve $x^d+y^d+z^d=0$ of degree
$d=4$ and $5$,
and $G$ consists of the automorphisms $(x, \ y, \ z)\rightarrow (\mu x,\ \nu
y,\ z)$,
where $\mu$, $\nu$ are  $\text{exp}({2\pi i\ov d})$.
(In this case, $g=(d-1)(d-2)/2$, $\vert G \vert =d^2$.)
\endproclaim
\demo{Proof} If $g(X)\ge 1$, by (4.2) we get the result.

 If $g(X)= 1$,  then the commutativity of $G$ implies $k\ge 2$.
By (4.2) we get $\tau \ge 1$, that is, $\vert G \vert\le 2g-2$.

Now we can assume   $g(X)=0$.
If $k\ge 6$,
by (4.2) we get $\tau \ge 1$, that is, $\vert G \vert\le 2g-2$.
Hence we can assume $k\le 5$.
We prove the case $k=3$;
for the case $k=4$ and $5$,
the proof is similar,
and left to the reader.

Clearly, by (4.2) we can assume $r_3\le 8$.
On the other hand, we claim $r_3 \ge 3$.
Otherwise,
let $Y=C/G(q_3)$.
Using Hurwitz formula to $Y \rightarrow X$,
which is unramified except $q_1 ,\ q_2 \in X$,
we have $Y \simeq \Bbb P^1$,
and $C$ is a hyperelliptic curve. A contradiction.

Let $$h_i={r_ir_3\ov\vert G \vert} \ ,\ \ \ \ \tau =1,\ 2.$$
By Lemma 4.3,
we have $h_i \in \Bbb N$,
and $h_1\ge h_2$
by our assumption.
By (4.2) we have
$$\tau =1-{r_3\ov h_1\vert G \vert}-{r_3\ov h_2\vert G \vert}-{1
\ov r_3}.\tag 4.5 $$
If $r_3=3$,
by (4.5) we have $\tau \ge {2\ov 3}-{6\ov\vert G \vert}$;
consequently, $\vert G \vert\le 3g+6$.

Now we can assume $4\le r_3 \le 8$.
By (4.2), we have  $\tau \ge 1-{3\ov r_3}$;
consequently, $$\vert G \vert\le {r_3\ov (r_3-3)} (2g-2).\tag 4.6$$

{\it Case 1. $h_2\ge 2$.}
By (4.5) we have $\tau \ge 1-{1\ov r_3}-{r_3\ov\vert G \vert}$;
consequently, $$\align \vert G \vert & \le {r_3\ov (r_3-1)} (2g-2+r_3)
\cr &  \le {4\ov 3}(2g-2+r_3) \cr & \le 3g+6\endalign$$
provided $g\ge 4r_3-26$.

Hence if either $r_3\le 7$ or $r_3=8$ and $g\ge 6$,
we get the result.
If  $r_3=8$ and $g < 6$,by (4.6) we have
$$\vert G \vert\le {16\ov 5}(g-1) <3g+6.$$
{\it Case 2. $h_2=1$, $h_1\ge 2$.}
Since $\vert G \vert=r_2r_3\ge r_3^2$ in this case,
we get $$g\ge {1\ov 2} (r_3^2-3r_3+2). \tag 4.7$$
By (4.5) we have $$\align \vert G \vert\le &
{r_3\ov (r_3-1)} (2g-2+{3r_3\ov 2})
 \cr <& 3g+6.\endalign$$
{\it Case 3. $h_1=h_2= 1$.}
By (4.5) we have $$\align \vert G \vert & = {r_3\ov (r_3-1)} (2g-2+2r_3)
\tag 4.8 \cr  & \le 3g+6\endalign$$
provided $g\ge 2(r_3-1)$.

By (4.7) we get the result when $r_3\ge 6$.

When $r_3=5$,
by (4.7) we have $g\ge 6$,
and  by (4.8) we have $g\equiv 0\pmod 2$.
Hence we have either $\vert G \vert\le 3g+6$
or $$(g,\ \vert G\vert ,\ r_1,\ r_2,\ r_3)=(6,\ 25 ,\ 5,\ 5,\ 5).$$
In the latter case,
since $g(X)=0$,
we can choose an isomorphism between $X$ and $\Bbb P^1$,
which maps $q_1$, $q_2$ and $q_3$ to the points $0$, $1$ and $\infty$
of $\Bbb P^1$.
Now it is easy to verify that the above $C$ is isomorphic
to the Fermat curve
 $x^5+y^5+z^5=0$,
and $G$ is the same as given in Theorem 4.4.

When $r_3=4$,
by (4.7) we have $g\ge 3$,
and  by (4.8) we have $g\equiv 0\pmod 3$.
Hence we have either $\vert G \vert\le 3g+6$
or $$(g,\ \vert G\vert ,\ r_1,\ r_2,\ r_3)=(3,\ 16 ,\ 4,\ 4,\ 4).$$
In the latter case, $C$ is isomorphic
to the Fermat curve of degree $4$.
\enddemo

{\bf Examples 4.9}. The estimates of Theorem 4.4 is best possible.
Let $m\ge 2$.
Let $C$ be given by the complete nonsingular model of affine equation $
y^3=x^{3m}-1$,
and $G$ consist of $x\rightarrow \mu x$;
 $y\rightarrow \nu y$,
where $\mu^{3m}=1$ and $\nu^3=1$.
Clearly, $C$ is nonhyperelliptic.
We have $g(C)=3m-2$, and $\vert G \vert =9m=3g(G)+6$.

We say a smooth curve $C$ is {\it bi-elliptic},
if $C$ can be represented as a ramified double covering of an elliptic curve.

\proclaim{Theorem 4.10}
Let $f\colon S \rightarrow B$ be a fibration of variable moduli.
Let $G$ be an abelian automorphism group of $S$,
inducing trivial action on $B$.
Then

(i) $\vert G\vert\le 4g-4$ ;

(ii) If the general fiber is neither hyperelliptic nor bi-elliptic,
then $\vert G\vert\le 3g-3$;

(iii) Suppose that $G$ is cyclic. Then  $\vert G\vert\le 2g+2$.\endproclaim
\demo{Proof}
Let $C$ be a general fiber of $f$.
We have a finite abelian covering
$$\pi \colon C \rightarrow X:=C/G .$$
Let $k$, $q_i$, $r_i$, $\tau$ and $G(q_i)$ be as in Notation 4.1.

 If $g(X)\ge 1$,
by (4.2) we get $\tau \ge 1$, that is, $\vert G \vert\le  2g-2$.

We can assume   $g(X)=0$. In this case we have $k\ge 4$
by the hypothesis of variable moduli.

When  $k\ge 6$,
by (4.2) we get $\tau \ge 1$, that is, $\vert G \vert\le  2g-2$.

When  $k=5$,
 we get either  $\tau \ge 2/3$, that is, $\vert G \vert \le 3g-3$
or  $\tau =1/2$, that is, $\vert G \vert = 4g-4$.
In the latter case,
$r_1= ...=r_5=2$.

Now we assume $k=4$.
Using Lemma 4.3 and  (4.2),
it is easy to see that if $r_3\ge 3$,
then  $\tau \ge 2/3$, that is, $\vert G \vert \le 3g-3$.

Hence we  assume $r_3=r_4=2$.
If $r_2\ge 6$, then
 $\tau \ge 2/3$, that is, $\vert G \vert \le 3g-3$.
If $r_2=5$,
by (4.2) and Lemma 4.3,
we have either  $\tau \ge 2/3$, that is, $\vert G \vert \le 3g-3$
or $(r_1,\ r_2,\ r_3,\ r_4)=(5,\ 5,\ 2,\ 2).$
In the latter case, we have  $\vert G \vert \ \vert\  10$, and
$g=4$ when $\vert G \vert=10$. Hence,  $\vert G \vert < 4g-4$.
If $r_2=4$,
by (4.2) and Lemma 4.3,
we have either  $\tau \ge 2/3$, that is, $\vert G \vert \le 3g-3$
or $$(r_1,\ r_2,\ r_3,\ r_4)\in \{(8,\ 4,\ 2,\ 2), \ (4,\ 4,\ 2,\ 2)\}.$$
In the latter case, we have $\tau \ge 1/2,$ that is,  $\vert G \vert \le 4g-4$.
If $r_2=3$,
by  Lemma 4.3 (i),
we have $\vert G \vert\  \vert\  12$.
In this case, we can assume
 $\vert G \vert = 12$.
Then  we have either  $\vert G \vert \le 3g-3$ (when $g \ge 5$)
or $(g,\ r_1)=(4,\ 6)$.
(Note that $(g,\ r_1)=(3,\ 3)$ is impossible by Lemma 4.3.)
If $r_2=2$,
by  Lemma 4.3 (i),
we have $\vert G \vert \ \vert\  8$.
Hence  we get  $\vert G \vert \le 3g-3$.
(Note that In this case, we have
 $g>3$ by (4.2).)

Summing up,
we have  $\vert G \vert \le 4g-4$; moreover,  $\vert G \vert \le 3g-3$
with the exceptional $\{g,\ \vert G\vert, \ k, \ (r_1, \ ... ,\ r_k) \}$ below:
$$\align  & \{5,\ 16,\ 5,\ (2,\ 2,\ 2,\ 2,\ 2)\},
\  \{4,\ 10,\ 4,\ (5,\ 5,\ 2,\ 2)\};\
\cr & \{6,\ 16,\ 4,\ (8,\ 4,\ 2,\ 2)\};\ \{4,\ 12,\ 4,\ (6,\ 3,\ 2,\ 2)\};\
\cr & \{3,\ 8,\ 4,\ (4,\ 4,\ 2,\ 2)\}.\endalign$$
Now we finish the proof of this theorem by proving the following claim.

{\it Claim.} If $\pi\colon
 C \rightarrow C/G$ is one of the above exceptional cases,
then $C$ is either hyperelliptic or bi-elliptic.

{\it Proof of the claim}. Let $Y=C/G(q_k).$
Then $\vert G(q_k) \vert=2$.
Using Hurwitz formula to $Y \rightarrow X$,
we have $g(Y)\le 1$.\qed

(iii) follows from Lemma 4.3 (iv).
See \cite{X3, Appendix A, Lemma A3} for a proof.
\enddemo

{\bf Remark.} It is likely that (i) and (ii) do not give the best bound when
$g\gg 0$,
but here we only need the result for small $g$.

 \heading \S 5. Abelian automorphism groups of small genus fibrations
\endheading

In \cite{Ch1},
Chen has shown that under the condition $K^2\ge 5$
and $S$ is not a product of two curves of genus $2$,
the order of an abelian subgroup of $\text{\text{Aut}(S)}$
is at most $12.5K^2+100$.
The purpose of this section is to give
a similar estimate for surfaces with small genus $3$, $4$,
and $5$.
As a consequence,
we give an estimate of the order of an abelian subgroup of $\text{Aut(S)}$
for the surface $S$ whose $1$-canonical map is composed with a pencil.

Let me start by giving a convenient lemma,
which plays an important role in this section.
For a proof, see \cite{X3, Appendix A}.

\proclaim{Lemma 5.1 (\cite{X3,  Lemma A1})}
 Let $f \colon S \rightarrow C$ be a fibration of
constant moduli,
of which the general fibers are of genus $g\ge 2$.
Let $F'$ be a singular fiber of $f$.
Denote by $\chi_{\text{top}}(F')$ the Euler?
topological character of $F'$.
Then
$$ \chi_{\text{top}}(F')+2g-2\ge \cases g-1, \ \   \text{if} \ g \text{
is odd}, \text{ and} \ F' \ \text{is a double} \cr\ \ \ \ \ \ \ \
\ \ \   \text{ curve of genus }
\ {(g+1)\over 2};\cr
g+2, \ \ \text{otherwise.}\endcases$$
In particular,
we have $\chi_{\text{top}}(F')+2g-2 \ge 4$ except if $g=3$,
and $F'$ is a double curve of genus $2$.\endproclaim

Let $S$ be a surface of general type with a pencil
$\Lambda$,
and $G$ an abelian group of automorphisms of $S$ whose elements
map $\Lambda$ onto itself.
Suppose $S$ has no ($-1$)-curves contained in a  fiber of $\Lambda$,
and let $g$ be the genus of a general element of $\Lambda$,
with $g\ge 2$.
Blowing up the base points (if necessary) of $\Lambda$,
$\Lambda$ corresponds to a relatively minimal fibration
 $f \colon S \rightarrow C$.
By hypothesis,
we have a homomorphism of $G$ into $\text{Aut(C)}$.
Let $H$ be its kernel.
There is an exact sequence
$$0 \rightarrow H \rightarrow G \rightarrow I\rightarrow 0, $$
where $I$ is the image of $G$ in  $\text{Aut(C)}$.
Let $F$ be a general fiber of $f$.
Then $H$ is an abelian  subgroup of  $\text{Aut(F)}$,
and $H$ induces trivial action on $C$.
We have  $ \# G=  \# H \#I$,
$ \# H \le  4g+4$,
and
$$ \# I\le  \cases 4b+4, \ \  \text{if} \ b\ge 2;\cr
\text{max}\{6,\ n \}, \ \ \text{if} \  b=1; \cr
\text{max}\{4,\ n \}, \ \ \text{if} \  b=0,\endcases \tag 5.2$$
where $b=g(C)$,
and $n$ is the number of singular fibers of $f$.

Let $\tau:=\text{min} \{\chi_{\text{top}}(F')+2g-2 \}$.
Then $\tau\ge 1$,  and
we have
$$ \align n\le  & {1\over \tau} \sum (\chi_{\text{top}}(F')+2g-2) \tag 5.3
\cr
\le  & {1\over \tau}[c_2(S)-4(g-1)(b-1)]\cr
\le  & {1\over \tau}({2g+1 \over g-1}K^2_S-16(g-1)(b-1)-24(b-1))\endalign $$
according to \cite{X4},
where the sum is taken over singular fibers $F'$ of $f$.

\proclaim{Proposition 5.4} Let $S$ be a minimal surface
of general type having two pencils $\Lambda_i$ ($i=1,\ 2$)
of genus $g\ge 2$.
Let $G$ be an abelian subgroup of $\text{Aut}(S)$.
Suppose that $K_S^2>4(g-1)^2$.
Then $\# G \le 16 K_S^2$.\endproclaim

\demo{ Proof} By assumption,
we have that $S\simeq C_1\times C_2$
for some smooth curve $C_i \in \Lambda_i$,
$i=1,\ 2$ (cf. \cite{X5, Proposition 6.4}).
Since the order of an abelian group of automorphisms of
genus $g$ is at most $4g+4$,
we have $\#G\le 32(g+1)^2$
(cf. \cite{X1, Example 2}).
Since $K_S^2=8(g-1)^2$ in this case,
we get the result.\qed\enddemo

\proclaim{Proposition 5.5} Let $S$, $G$ be as above, and that
 $S$ has a  pencil
$\Lambda$ of curves of genus $3$.
Suppose that $K_S^2>16$,
and $K_S^2\ge 3\chi(\Cal O_S)$.
Then
$$\# G\le 24K^2+64.$$\endproclaim
\demo{ Proof} We note that $\Lambda$ has no base point:
we have $(mK-nC)C=0$ for $C \in \Lambda$,
and  some positive integral numbers $m$, $n$.
Then the Hodge's index theorem implies
$K^2\le 9$,
a contradiction.

Let  $f \colon S \rightarrow C$ be the fibration of genus $3$ associated
with $\Lambda$, and $b=g(C)$.
Since $K^2>16$, by Proposition 5.4,
we can assume that $S$ has only one fibration of genus $3$.
Hence each element of $G$ maps $\Lambda$ onto itself.

If $f$ is of variable moduli,
then by lemma 2 we have
$\# H\le 8$.
By (5.2), (5.3) and the assumption we have
$\#I\le 3K_S^2+8$,
and consequently $$\# G\le 24K^2+64.$$

If $f$ is of constant  moduli,
Lemma 5.1 allows us to conclude $\# I\le {1\over 5}(c_2(S)-4(g-1)(b-1) )$;
consequently, $$\# G \le {56\ov 5}( K^2+ 16) $$
except when there are many double curves of genus $2$.
In the latter case let $D$ be a general divisor in the moving part of
$\vert K \vert$.
Then $0\le DF\le4$.

Now we consider two cases separately.

{\it Case 1. $DF>0$.}
Consider the projection $\psi  \colon D \rightarrow C$
induced by $f$.
$\psi$ must be ramified along each double fiber.
Let $m$ be the number of double fibers of $f$.
By Hurwitz formula,
we have
$$2K^2\ge 2p_a(D)-2\ge (DF) (2b-2+{m \over 2});$$
so $m\le 2K^2+4$,
and by Lemma 5.1, (5.2) and (5.3) we have
$$\# G \le (4g+4){m \over 2}\le 16K^2+32.$$

{\it Case 2. $DF=0$.}
In this case,
$\vert K \vert$ is composed with a pencil of genus $3$,
and $f=\varphi _K$, the $1$-canonical map of $S$.
Let $\vert K \vert=\vert M \vert+Z$,
where $Z$ is a fixed part of $\vert K \vert$,
and $\vert M \vert$ is the moving part of $\vert K \vert$.
Then $M\sim aF$,
where $a\ge p_g-1$ (see \cite{Be, proposition 2.1}),
and $$K^2\ge aKF\ge (p_g-1)KF \ge 4\chi -8.$$
Hence we have
$$c_2(S)-4(g-1)(b-1)=12\chi-K^2-4(g-1)(b-1)\le 2K^2+32; $$
consequently, by Lemma 5.1 and (5.3),
$\# I \le K^2+16, $
and  we get  $ \# G \le 16(K^2+16)$.\qed
\enddemo

\proclaim{Proposition 5.6} Notations being as in proposition 5.5.
Suppose  $K^2>36$, and that $S$ has a pencil
 $\Lambda$ of curves of genus $4$ .
Suppose that $K_S^2>36$,
and $K_S^2\ge 4\chi(\Cal O_S)$.
Then
$$\# G\le 24K^2+144.$$\endproclaim
\demo{ Proof}It is easy to see  that $\Lambda$ has no base point
as in the proof of proposition 5.5.

Let  $f \colon S \rightarrow C$ be the fibration of genus $4$ associated
with $\Lambda$, and $b=g(C)$.
Since $K^2>36$,  by Proposition 5.4,
we can assume that $S$ has only one fibration of genus $4$.
Hence each element of $G$ maps $\Lambda$ onto itself.

If $f$ is of constant  moduli,
then by Lemma 5.1  we have
 $$\# G\le 10K^2+240.$$

If $f$ is of variable  moduli,
Lemma 4.10 allows us to conclude $\# H\le 12$;
consequently, $\# G \le 24K^2+144 $
as in the proof of Proposition 5.5.\qed
\enddemo

\proclaim{Proposition 5.7} Notations being as in proposition 5.5.
Suppose that $K_S^2>64$,
 $K_S^2\ge {24\ov 5}\chi(\Cal O_S)$, and that $S$ has a pencil
 $\Lambda$ of curves of genus $5$.
Then
$$\# G\le 24K^2+256.$$\endproclaim
\demo{ Proof} The proof is similar to that of Proposition 5.6,
and left to the reader.\qed
\enddemo

\proclaim{Proposition 5.8}  Let $S$ be a minimal smooth surface of general
type,
and $K$ the canonical divisor of $S$.
Let $G$ be an abelian group of automorphisms of $S$ (i.e. $G\subset
\text{Aut} (S)$).
And let $\varphi_1$ be the $1$-canonical map of $S$.
Suppose that $\chi(\Cal O_S)\ge 21$,
and that  $\text{dim}\ \varphi_1(S)=1$.
Then $$\# G \le 12.5K^2+469.$$
\endproclaim
\demo{Proof}
By assumption,
$\varphi$ is composed with a pencil of curves of genus $g$
with $2\le g\le 5$,
and the pencil has no base points according to \cite{Be}.
Let $f\colon S\r C$ be a fibration of genus $g$ induced by  $\varphi_1$.
Let $F$ be a general fiber of $f$.
If $g=2$, then $\#G\le12.5K_S^2+100$ \cite{Ch}.
Hence we can assume that $g\ge 3$.

Let $\vert K_S\vert=\vert M \vert+Z $,
where $Z$ is the fixed part of $\vert K_S \vert$,
and $\vert M \vert$ is the moving part of $\vert K_S \vert$.
Then $M\sim aF$,
where $a\ge p_g-1$ (see \cite{Be, Proposition 2.1}),
and $$K_S^2\ge aK_S.F\ge (2g-2)(p_g-1)\ge (2g-2)(\chi-2).$$
Consequently, $\#I\le {1\ov \tau } (({6\ov g-1}-1)K_S^2+24+4(g-1))$,
where $\tau$ is as in (5.3).
Now it is easy to see that
$$\#G\le \cases  8K^2+640\le  12K^2+496, \ \ \text{if}\ \  g=4;\cr
 12K^2+432, \ \  \text{if}\ \  g=5.\cr\endcases$$
When $g=3$, one has $ K_S^2\ge {21\ov 4}p_g-{71\ov 6}$ \cite{Su}.
Consequently, $\#I\le {1\ov \tau } (({9\ov 7}K_S^2+47+{1\ov 21})$,
 and $\#G\le {72\ov 7}K^2+376+{8\ov 21}$.\qed \enddemo

\heading \S 6. Abelian subgroups for surfaces with large $K^2_S$ \endheading

Our main result  in  this section
is  the  following.

 \proclaim{Theorem 6.1} Let $S$ be a minimal smooth surface of general
type over the complex number field,
and $K$ the canonical divisor of $S$.
Let $G$ be an abelian group of automorphisms of $S$ (i.e. $G\subset
\text{Aut} (S)$). Suppose that the  dimension of
 the $1$-canonical image of $S$ is
 $2$, $\chi \ge 14$, $K^2 \ge 82$,
and   that $S$ has no pencil of curves of genus $\le 5$.
Then $$\# G \le 24K^2 +16.$$
Moreover,
if  the $1$-canonical map  of $S$ is birational,
then  $$\# G \le 18K^2 +18.$$  \endproclaim

The arguments here are inspired by the work  of Xiao  \cite{X1}.
We consider the natural action of the abelian group $G$
on the space $H_n=H^0(S,\ nK_S)$,
for a fixed positive integer $n$.
Because $G$ is finite abelian,
such an action is diagonalisable,
in other words $H_n$ has a basis composed of semi-invariant vectors.
Consider two such semi-invariants  $v_1,\ \ v_2$ in $H_n$,
with
$$\sigma (v_i)=\alpha _i(\sigma)v_i \ \ \text{for}\ \ \sigma \in G,$$
where
$ \alpha_i$ are the corresponding characters of $G$.
Suppose  that the  two  semi-invariants  $v_1,\ \ v_2$
corresponding to a same character of $G$ (i.e., $\alpha_1=\alpha_2$ ),
 and let
 $D_1$ and $D_2$ be  the corresponding
divisors  in $ \vert nK_S\vert$.
Then
 $D_1$ and $D_2$ generate a pencil $\Lambda$
whose general fiber $F$ is fixed by $G$.
Therefore $\#G$  is limited by the order of the group of automorphisms of the
normalization $\tilde F $ of $F$ as a smooth curve.
But $\#Aut(\tilde F)$ increases proportionally with $K^2_S$,
as $g(\tilde F)$ so does.

Instead of considering the  natural map
$$ H_n  \otimes H_n    \rightarrow H_{2n}, $$
We may consider the natural map
$$ H_{n-t}  \otimes H_{n+t}  \oplus H_n  \otimes H_n  \rightarrow H_{2n}, $$
which is compatible with the above actions of $G$:
i.e.,  if   $v_1 \in H_{n-t}$, $v_2 \in H_{n+t}$ (resp. $w_i \in H_n$
$i=1,  \ 2$)
are  two  semi-invariants,
then $v_1 \otimes v_2$ (resp.  $w_1 \otimes w_2$)
is  semi-invariant in $H_{2n}$.
If there are more than $\text{dim}\  H_i$ semi-invariants in
$H_i$ for some $i\le {n+t}$,
then  there are semi-invariants in
$H_i$ (therefore in $H_{2n}$) with the  same character,
and we are done.
So we may assume that  there are exactly  $\delta_i=\text{dim}\  H_i$
($i\le {n+t}$) semi-invariants $v_j^i$ ($j=1,\ ...,\ \delta_i$) in
$H_i$,
corresponding to mutually different characters.
Each vector $v_j^i$
corresponds to a unique divisor $D_j^i$ in $\vert iK \vert$.
The relation $v_j^i \otimes v_k^l=c v_m^{r}\otimes v_n^{s}$
(where $c$ is a constant,
and $i+l=r+s=2n$) in $H_{2n}$ translates to a  relation
$$D_j^i+ D_k^l= D_m^{r}+D_n^{s} \tag *$$
between these divisors.

Fix an semi-invariant $u\in H_t$ when $H_t\not= 0$.
$u$ corresponds to a unique divisor $U$ in $\vert tK \vert$ which is fixed
by $G$.
Then we can consider the finite set  $\Sigma _i$
  of points corresponding to $D_j^i$
in a certain divisorial space $P_i$ defined in  \cite{X1, $\S 1$}
(see also \S 3),
and there are  natural embeddings:
$$ \Sigma _{n-t} \rightarrow  \Sigma _n \rightarrow  \Sigma _{n+t},  \ \ \ \
P_{n-t} \overset {l_2} \to  \rightarrow
P_n \overset {l_3} \to  \rightarrow P_{n+t}$$
defined by $U$.
In such a setting,
an semi-invariant in $H_{2n}$ of the form  $v_j^i\otimes v_k^l$
corresponds naturally to the mid-point of two points in
$ \Sigma _{n+t}$
corresponding to $D_j^i$ and $ D_k^l$,
and a relation of the form (*) means that the corresponding mid-points
coincide.

Denote by $ \Cal S_1$ the set of
mid-points of two points $p$, $q$ in $ \Sigma_ {n+t}$
such that either $p$ is in $ \Sigma_{n-t}$ and $q$ is in $ \Sigma_{n+t}$
 or $p$ and $q$
are in $ \Sigma_n$.
 We reduce the problem to that of comparing the number of points in $ \Cal S_1$
and  the dimension of $H_{2n}$.
In this way we show that for $n=2$ and  $t=1$
 the number of points in $ \Cal S_1$
is larger than   the dimension of $H_4$, provided
 that the  dimension of
 the $1$-canonical image of $S$ is
 $2$, $\chi \ge 14$, $K^2 \ge 82$,
and   that $S$ has no pencil of curves of genus $\le 5$.

Assume $H_1 \not=0$.
Fix an semi-invariant $u \in H_1$ for the natural action of $G$ on $H_1$.
Let $U$ be the divisor in $\vert K \vert$ corresponding to $u$.
We have natural maps:
$$H_1 \overset { \otimes  u} \to  \rightarrow
H_2 \overset { \otimes  u} \to  \rightarrow H_3;$$
$$ \vert  K \vert  \overset {+U} \to  \rightarrow
 \vert 2K \vert  \overset {+U} \to  \rightarrow  \vert  3K \vert ;$$
$$P'_1 \overset {l'_1} \to  \rightarrow
P'_2 \overset {l'_2} \to  \rightarrow P'_3.$$
If we take $[ nU ]$ be the origin of $P'_n$,
then $l'_i (i=1, \ 2)$ are embeddings of linear spaces.
Fix $P'_1$ and $P'_2$ as subspaces of $P'_3$ in this way in the future.
Since $U$ is fixed by $G$,
 $l'_i$ ($i=1,\ 2$) induce  natural embeddings:
$$ \Sigma _1 \rightarrow  \Sigma _2 \rightarrow  \Sigma _3,  $$
 and natural embeddings:
$$P_1 \overset {l_1} \to  \rightarrow
P_2 \overset {l_2} \to  \rightarrow P_3.$$
Let $v $ (resp.$ w)$ be an semi-invariant in $H_1$ (resp. $H_3$), and
 $p$ (resp. $q$) the corresponding
point in $ \Sigma_1$ (resp. $ \Sigma_3$).
Let $D$ be the divisor in $\vert 4K\vert$
corresponding to the vector $v \otimes w$.
Then $l_2 ([ {1 \over2}D ])$  corresponds to a point in $P_2$, which is just
the mid-point ${1 \over2}(l_2l_1(p)+q)$:
in fact,    Let $D'$ (resp. $D''$)
be the divisor in  $\vert K\vert $ (resp. $\vert 3K\vert$)
corresponding to the vector $v $ (resp. $w$).
Then $D=D'+D''$, and
$$\align  l_2([ {1 \over2}D ]) & =[ {1 \over2}D+U ]
\cr & =[ {1 \over2}(D'+2U)+{1 \over2}D'' ]  ={1 \over2}[ ( D'+2U)+D''-D_1 ]
\cr & ={1 \over2}([ D'+2U ] +[ D'' ] )  ={1 \over2}(l_2l_1(p)+q)
.\endalign $$
Similarly, let $v_1,\ \ v_2$ be two  semi-invariants in $H_2$,
 and
 $p$ , $q$ the corresponding
points in $ \Sigma_2$.
Let $D$ be the divisor in $\vert 4K\vert$
corresponding to the vector $v_1 \otimes v_2$.
Then $l_2 ([ {1 \over2}D ])$  corresponds to a point in $P_2$, which is just
the mid-point ${1 \over2}(l_2(p)+l_2(q))$.

{\bf  Remark.} $\Sigma _n$ is  integrally convex
with respect to the lattice  $\Cal  L$
consisting of points corresponding to divisors linearly equivalent to
$nK$  \cite{X1, Lemma 3}.
Clearly it is easy to verify that
$\Sigma _i $ is relatively convex in $\Sigma _{i+1} $
for all $i$,
as we consider $\Sigma _i $ as a subset of $\Sigma _{i+1} $
in the above way.

\proclaim{Lemma 6.2 }
Suppose  $K^2  \ge 82$,
and let $\Sigma_i$ be a basic set in $H_i$ for $i=2$, $3$.
Assume that either there is a chain in $\Sigma _3$ of length
$ \ge  {1 \over 10}  \# \Sigma_3$ or there is a
chain in $ \Sigma_2$ of length $  \ge  {1 \over 5}  \# \Sigma_2$.
Then $S$ has a pencil of curves of genus $\le 5$.
 \endproclaim
\demo{Proof} The proof is a modification of \cite{X1, lemma 4}.
Suppose there is a chain in   $\Sigma _3$ of length
$ \ge  {1 \over 10}  \# \Sigma_3$.
(For the latter case, the proof is similar.)
Let $k$ be the length of the  longest chain  in $\Sigma _3$.
Then there is a pencil $\Lambda$ of curves on $S$ such that there is a divisor
in $\vert 3K \vert$ containing $k$ times a fiber $F$ of $\Lambda$.
We have $KF\le 9$.
Since $K^2\ge 82$,
the Hodge index theorem implies $F^2=0$.
Now the evenness of $F^2+KF$ forces $KF\le 8$,
which means that $\Lambda$ is a pencil of curves of genus $\le 5$.
\qed \enddemo

 \proclaim{Proposition 6.3} Let $S$ be a minimal surface whose $1$-canonical
image is of dimensional $2$.
 Suppose  $\chi \ge 14$
and that $S$ has no pencils of curves of genus $\le 5$.
Then $H_4$ is not uniquely decomposable.
  \endproclaim
\demo{Proof}We can suppose $H_n$ is  uniquely decomposable for $n\le 3$,
for otherwise the corollary is trivially true.
Let $\Sigma_i$ be a basic set in $P_i$ ( $i=1, \ 2,\ 3$).
By \cite{X1, lemma 3} and by the hypothesis that
 $1$-canonical
image is of dimensional $2$,
we have the dimension of $\Sigma _1$ is $\ge 2$;
also by  \cite{X1, lemma 3},
we have the dimension of $\Sigma _2$ is $\ge 3$;
By lemma 6.2,
the condition of  lemma 2.7 is satisfied for
$$ \Sigma _1 \subset  \Sigma _2 \subset  \Sigma _3.  $$
Then this corollary  results from the relations
$$ \align \#\Sigma_i= & \text{dim}\ H_i
={i(i-1)\over 2}K^2+\chi \ \ \ (i=2,\ 3),\cr
K^2\le &  9\chi \ \ \ (\text{ Bogomolov- Miyaoka-Yau's
inequality}),\endalign$$
that $\# ( \Sigma _1 .  \Sigma _3 \cup
    \Sigma _2 .  \Sigma _2) >\text{dim}\  H_4.  $
In particular,
there are more than $\text{dim}\  H_4$ semi-invariants in $H_4$.
\qed \enddemo

 \proclaim{Lemma 6.4 } Suppose  $K^2  \ge 10$
and with  no pencils of curves of genus $2$.
Let $G$ be an abelian group of automorphisms of $S$,
such that $\vert 4K\vert$ contains a pencil $\Lambda  $
whose general numbers are fixed under the action of $G$.
Then $$\# G \le 24K^2 +16.$$
Moreover,
if  the $1$-canonical map  of $S$ is birational,
then  $$\# G \le 18K^2 +18.$$
\endproclaim
\demo{Proof}
Blowing up the base points of $\Lambda$,
we get a surface $S'$ such that  $\Lambda$ is associated to a fibration
$f\colon S'\rightarrow C$.
Let $F$ be a general fiber of $f$,
and let $k$ be the number of fibers of $f$ contained in a general
member of the moving part of  $\Lambda$.
We have
$$ g(F)-1\le {2k+8\over k^2}K^2.$$
Let $H$ be the stabilizer of $F$.
Then the index of $H$ in $G$ is at most $k$.
Because $H$ is an abelian group of automorphisms of curve $F$,
we have
$$ \# H\le 4g(F)+4.$$
Therefore if $k \ge 2$,
we get $$\# G \le \cases 36K^2 +24;\cr
 18K^2 +18, \ \ \text{if } \ f\ \text{is non-hyperelliptic.}\cr
\endcases $$
This allows us to assume $k\le 1$,
consequently there is no divisor in $\vert 2K\vert$ whose pull-back
on $S'$ contains $F$.

Let $\pi\colon F\rightarrow B:=F/H$ be the projection map.
Let $O_1$, ..., $O_l$ be the orbits of the action of $H$ on $F$
contained in pull-backs of fixed divisors in  $\vert 2K\vert$,
$n_i$ the number of points in $O_i$.
If $l<4$,
then the image of the bicanonical map of $S$ is a rational surface
(see e.g. \cite{X1, p 624} for a proof);
therefore  $S$ has a pencil of curves of genus $2$ \cite{1?},
contradiction to the hypothesis.
So we can assume $l\ge 4$.
Using Hurwitz formula to  $\pi$,
we get
$$ (l-2)\# H\le 2g(F)-2+\sum^l_{i=1} n_i.$$
Because $n_i\le 2KF \le 8K^2 / k$,
we see that $\#G \le 17K^2$,
once $l\ge 6$.
When $4\le l\le 5$,
remark that the pull-back $E$ of each fixed divisor
in  $\vert 2K\vert$ has a unique representation
$$ E=a_1O_1+\dots +a_lO_l$$
with $0\le a_in_i\le 8K^2$,
$i=1$,  $\dots$,
$l$.
Let $b_i$ ($i=1$, $\dots$,
$l$)
be the maximum of the $i$-th coefficient in all such  representations,
we must have
$$b_1.  \dots.   b_l \ge \text{dim}H_2=K^2+\chi,\ \
b_in_i\le 8K^2.$$
Now if $l=5$,
then $\sum {1\ov b_i}<4+1/K^2$;
consequently,
$$\sum n_i<8(4+{1\ov K^2})K^2, \ \ \text{and}\ \ \#G<{52\ov
3}K^2+{8\ov 3}.$$
When $l=4$,
taking also into account the fact that no such pull-back divisor
can contain another,
it is easy to see that the minimum for $\sum_{i=1}^4 n_i$ is
achieved when,
say $$\align & b_1=b_2={1\ov 2}\o{dim}H_2-1,\ \ b_3=b_4=1,\cr
& n_1=n_2=16, \ \ \ \ \ \ \ \ \ n_3=n_4=8K^2,\endalign$$
and consequently $\#G\le 18K^2+16$.
\qed \enddemo

\demo{Proof of theorem 6.1} Proposition 6.3 guarantees that there is a pencil
  $\Lambda$
in $\vert 4K \vert$ whose every element is a fixed divisor by $G$.
Hence by Lemma 6.4, we get the result.
\qed \enddemo

 \proclaim{Theorem 6.5} Let $S$ be a minimal smooth surface of general
type over the complex number field,
and $K$ the canonical divisor of $S$.
Let $G$ be an abelian group of automorphisms of $S$ (i.e. $G\subset
\text{Aut} (S)$). Suppose that
 $K^2 \ge 181$,
and   that either $S$ has no pencil of curves of genus $g$, $3\le
g\le  5$ or $K_S^2\ge {12(g-1)\ov g+5}\chi(\Cal O_S)$
when  $S$ has a pencil of curves of genus $g$,
$3\le g\le  5$.
Then $$\# G \le 24K^2 +256.$$
  \endproclaim
\demo{Proof} By a theorem of Beauville (cf.\cite{Be},
\cite{BPV}), we have that either $\varphi_1$ is composed with a
pencil of curves of genus $g$, $2\le g\le 5$,
and the pencil has no base points, or $\o{dim}\varphi_1(S)=2$.
Hence the result follows from Theorem 6.1 and Propositions 5.5,
5.6, and
5.7.\qed \enddemo

\heading \S 7. Abelian subgroups for small numerical invariants
 \endheading

 \proclaim{Theorem 7.1} Let $S$ be a minimal smooth surface of general
type over the complex number field,
 $K$ the canonical divisor of $S$,
and $\chi(\Cal O_S)$ the Euler characteristic of the structure sheaf of $S$.
Let $G$ be an abelian group of automorphisms of $S$ (i.e. $G\subset
\text{Aut} (S)$).
Then $$\# G \le 36K^2 +24,$$
provided $\chi(\Cal O_S)\ge 8$.  \endproclaim

\proclaim{ Lemma 7.2}
Suppose  $\chi \ge 5$, and that $P_6$ is uniquely decomposable.
 Then the dimension of $ \Sigma_3$ is  $3$.
\endproclaim
\demo{ Proof} Let $d$ be the dimension of $\Sigma_3$,
Then
$$ \# ( \Sigma_3.\Sigma_3)\ge (d+1)(\#\Sigma_3-{d\over2}) \ \ \  \text{
(cf. \cite{5, Lemma 1}}). $$
Now assume $d\ge 4$,
then taking into account the inequality $\#\Sigma_3\ge d+1$,
we have
$$ \#(\Sigma_3.\Sigma_3)\ge 5 \# \Sigma_3 -10.$$
On the other hand,
we have
$$\text{dim}\  H_6=15K^2+\chi,\ \  \text{dim}\  H_3=3K^2+\chi.  $$
So we get $\#(\Sigma_3.\Sigma_3)\ge 5\text{dim}\  H_3-10>\text{dim}\  H_6,$
when $\chi \ge 3$, a contradiction.
Hence we have $d\le 3$.
 Now $d\le2$ is impossible:
otherwise, the image of the $3$-canonical map of $S$ is either a
rational curve or a rational surface (cf. \cite{5, Lemma 3}).
 \qed \enddemo
 we may copy the proof of  \cite{X1, Lemma 4} to get the following.
\proclaim{Lemma 7.3}
Suppose  $K^2  \ge 10$,
and let $\Sigma_i$ be a basic set in $H_i$ for $i=3$, $4$.
Assume that either there is a chain in $\Sigma _4$ of length
$ \ge  {1 \over 6}  \# \Sigma_4$ or there is a
chain in $ \Sigma_3$ of length $  \ge  {1 \over 4}  \# \Sigma_3$.
Then $S$ has a pencil of curves of genus $2$.
\qed \endproclaim

 \proclaim{Lemma 7.4} Suppose  $\chi \ge 8$
and that $S$ has no pencils of curves of genus $2$.
Then $H_6$ is not uniquely decomposable.
  \endproclaim
\demo{Proof}We can suppose $H_n$ is  uniquely decomposable for $n\le 4$,
for otherwise the corollary is trivially true.
Let $\Sigma_i$ be a basic set in $P_i$ ( $i=2,\ 3,\ 4$).
By Lemma 7.2 and Lemma 7.3,
the condition of  Lemma 2.?5 is satisfied for
$$ \Sigma _2 \subset  \Sigma _3 \subset  \Sigma _4.  $$
Then this corollary  results from the relations
$$ \align \#\Sigma_i= & \text{dim}\ H_i
={i(i-1)\over 2}K^2+\chi \ \ \ (i=2,\ 3,\ 4),\cr
K^2\le &  9\chi \ \ \ (\text{ Bogomolov- Miyaoka-Yau's
inequality}),\endalign$$
that $\# ( \Sigma _2 .  \Sigma _4 \cup
    \Sigma _3 .  \Sigma _3) >\text{dim}\  H_6.  $
In particular,
there are more than $\text{dim}\  H_6$ semi-invariants in $H_6$.
\qed \enddemo

 \proclaim{Lemma 7.5 } Suppose  $K^2  \ge 10$
and with  no pencils of curves of genus $2$.
Let $G$ be an abelian group of automorphisms of $S$,
such that $\vert 6K\vert$ contains a pencil $\Lambda  $
whose general numbers are fixed under the action of $G$.
Then $$\# G \le 36K^2 +24.$$
\endproclaim

\demo{Proof}
We modify the proof of \cite{X1, lemma 7} to get our results:
we can assume that there is no divisor in $\vert 2K \vert$
 whose pull-back contains a general fiber $F$
of the fibration $f \colon S' \rightarrow C$
induced by a pencil $\Lambda$
of $G$-invariant divisors in  $\vert 4K \vert$.
Let $O_1$, ..., $O_l$  be the orbits on $F$
contained in pull-backs of fixed divisors  in the moving parts
of the pull-backs of    $\vert 2K\vert$.
As $K^2\ge 10$,
the $2$-canonical map  is birational
by Bombieri's theorem,
and so $l\ge 4$.
Using Hurwitz formula to  $\pi$,
we get
$$ (l-2)\# H\le 2g(F)-2+\sum^l_{i=1} n_i.$$
Because $n_i\le 2KF \le 12K^2 / k$,
we see that $\#G \le 34K^2$,
once $l\ge 5$,
and  $\# G\le 33K^2+24$
when $l=4$ as in \cite{X1, Lemma 7}.
\qed \enddemo

\demo{Proof of theorem 1} If $S$ has  a
relatively minimal genus $2$ fibration,
then $\# G \le 12.5K^2 +100$ (cf. \cite{Ch1, Theorem 0.2}).
Hence  we can suppose that $S$ has
no pencil of curves of genus $2$,
then  Lemma 7.4 guarantees that there is a pencil
  $\Lambda$
in $\vert 6K \vert$ each of  whose  elements is a fixed divisor by $G$.
Hence by Lemma 7.5, we get the result.
\qed \enddemo

As a consequence of Theorem 7.1,
we have that the order of an Abelian subgroup of $\text{Aut}(S)$
is at most  $\# G \le 36K^2 +24,$
provided $K^2\ge 64$.
Here we give a similar estimates for surfaces with $K^2<64$.

\proclaim{ Lemma 7.6}
 Suppose  $K^2  \ge 4$ (resp. $K^2  \ge 2$ ).
Then $H_{12}$ (resp. $H_{16}$) is not uniquely decomposable.
\endproclaim
\demo{Proof}
 We may consider the natural map
$$ H_4  \otimes H_8  \oplus H_6  \otimes H_6  \rightarrow H_{12} $$
$$ (\text{resp.}\  H_5
\otimes H_{11}  \oplus H_8  \otimes H_8  \rightarrow H_{16}) $$
instead of
$$ H_2  \otimes H_4  \oplus H_3  \otimes H_3  \rightarrow H_6. $$

  Note that if there is a chain in
$\Sigma_6$ (resp. $\Sigma_8$,
$\Sigma_{11}$) of length $\ge {1\over 4}\# \Sigma_6$
(resp.  $\ge {1\over 4}\# \Sigma_8$,  $\ge {1\over 6}\# \Sigma_{11}$)
then $S$ has a pencil of curves of genus $1$,
contradicting that $S$ is of general type (cf. \cite{X1, Lemma 4} for a proof).
We choose and fix an semi-invariant $u'$ in $H_2$ (resp. $u''$ in $H_3$)
instead of $u$ in $H_1$, as in $\S 1$.
Then one checks immediately that the first half of $\S 6$ goes
 for the new pair.
\qed \enddemo
 \proclaim{Theorem 7.7} Let $S$ be a minimal smooth surface of general
type over the complex number field,
and $K$ the canonical divisor of $S$.
Let $G$ be an abelian group of automorphisms of $S$ (i.e. $G\subset
\text{Aut} (S)$).
Then $$ \# G
\le  \cases
 114K^2 +24, \ \text{provided }\ \ 4 \le K^2\le 63;\cr
 200K^2 +22, \ \text{provided }\ \ 2 \le K^2\le 3;\cr
270, \ \text{provided }\ \  K^2=1. \endcases$$
  \endproclaim
\demo{Proof}
If  $K^2\ge 4$ (resp. $K^2\ge 2$),
by Lemma 7.6, $H_{12}$ (resp $H_{16}$) is not uniquely decomposable.
We modify the proof of \cite{X1, Lemma 7} to get our results:
we can assume that there is no divisor in $\vert 3K \vert$
(resp.  $\vert 4K \vert$) whose pull-back contains a general fiber $F$
of the fibration $f \colon S' \rightarrow C$
induced by a pencil $\Lambda$
of $G$-invariant divisors in  $\vert 12K \vert$ (resp. $\vert 16K \vert$).
Let $O_1$, ..., $O_l$  be the orbits on $F$
contained in pull-backs of fixed divisors  in the moving parts
of the pull-backs of    $\vert 3K\vert$ (resp. $\vert 4K \vert$).
As $K^2\ge 4$ (resp. $K^2\ge 2$),
the tricanonical map (resp. $4$-canonical map) is birational
by Bombieri's theorem,
and so $l\ge 4$.
And as
$$  g(F)-1\le {6k+72\over k^2}K^2$$
$$(\text{resp.}\  g(F)-1\le {8k+128\over k^2}K^2),$$
 $n_i\le 36K^2 / k$ (resp. $\le 64K^2 / k$),
we get $\#G \le 112K^2$ (resp. $\le 198K^2 $),
once $l\ge 5$,
and $\#G \le 114K^2+24$ (resp. $\le 200K^2+22$)  when $l=4$
as in the proof of Lemma 7.5.

If $K^2=1$, the proof of \cite{X2, Appendix A, Theorem A1}
shows that $\#G\le 270$.
\qed \enddemo

\heading \S 8. Miscellaneous results
 \endheading

 \proclaim{Theorem 8.1} Let $S$ be a simply connected even
 smooth surface of general
type over the complex number field, i.e., the intersection form
on $H^2(S,\ \Bbb Z)$ is even,
or equivalently,
$K_S=2L$ for some integral divisor $L\in \o{Div}(S)$.
Let $G$ be an abelian group of automorphisms of $S$ (i.e. $G\subset
\text{Aut} (S)$).
Suppose that

(i)$ \phi_L$ is generally finite,
and $ \phi_{2L}$ is birational;

(ii) $S$ has no fibrations of genus $3\le g\le 8.$

(iii)$ K_S^2>196.$ Then
 $$\# G \le 12K^2 +24.$$
  \endproclaim
\demo{Proof} Since $S$ is simply connected,
for any $\sigma \in G$,
we have $\sigma^*L\equiv L$.
Therefore there is a natural action of the abelian group $G$ on
$H_n(L):=H^0(S,\ nL)$,
for a positive integer $n$.
Hence we can use $nL$ (instead of $nK$) to define the space $P_n(L)$
and the basic set $\Sigma_n(L)$ as in \S 3.

Then we can modify the proof of Theorem 6.1 to get the result.
We leave the detail of the proof to the reader.
\qed \enddemo

 \proclaim{Corollary 8.2} Let $d_i \in \Bbb Z$,
$d_i\ge 2$,
for $i=1,\ \cdots,\ N-2$,
  let $S$ be a
 smooth complete intersection of type $(d_1,\ \cdots,\ d_{N-2})$
in $\Bbb C\Bbb P^N$, and
let $G$ be an abelian group of automorphisms of $S$ (i.e. $G\subset
\text{Aut} (S)$).
Suppose that
$\sum_{i=1}^{N-2}d_i-N$ is odd,
and $ K_S^2>196.$ Then
 $$\# G \le 12K^2 +24.$$
  \endproclaim

\demo{Proof} By Theorem 8.1,
it is enough to show that $S$ has no fibrations of genus $\le 8$.
Otherwise,
suppose that $S$ has a fibration of genus $g \le 8$.
Let $F$ be a general fiber of such a fibration.
Since $\o{deg}F\ge N$ (cf. \cite{Mu, p 77}),
we have $$14 \ge 2g-2\ge N(\sum d_i-N-1)\ge N(N-5).$$
Hence $N\le 7$, and $\sum d_i\le N+1+[{14\ov N}]$.
Now we can verify case-by-case that this is impossible since
$ K_S^2>196.$\qed\enddemo

For surfaces in  $\Bbb C\Bbb P^3$,
we have a better estimation for the order of abelian
automorphism groups.

 \proclaim{Theorem 8.3} Let $G$ be an  abelian group
 of automorphisms of a  surface $S$ in  $\Bbb C\Bbb P^3$
of degree $d\ge 5$.
Then $\#G\le 3d^2(d-2)+9.$\endproclaim

\demo{Proof} Changing coordinates,
we can assume that each element $\sigma\in G$
has the form $\o{diag}(*,*,*,*)$
since the finite  abelian group $G$ can be diagonalisable.
Let $a_1m_1(x)+\dots +a_tm_t(x)$ be the defining equation
of $S$,
where $m_i(x)$ are monomials of degree $d$ in $\Bbb C[x_0,\
\dots,\ x_3]$,
and $a_i\in \Bbb C$.
Since $S$ is smooth,
there exists a surface $Y$,
 the defining equation given by a general linear combination of
$m_i(x)$,
such that $Y$ and the complete intersection $C:=X\cap Y$
are both smooth.
We have a natural group map $\beta\colon G\r \o{Aut}(C)$.
$\o{Ker}\beta=1$ since $C$ is not contained in a hyperplane of
$\Bbb P^3$.
Since $C$ is non-hyperelliptic,
we have the order of an abelian automorphism group is $\le 3g(C)+6$
(Theorem 4.9),
where $g(C)=d^2(d-2)+1$.
Consequently, we have $\#G\le 3d^2(d-2)+9.$
\qed \enddemo

It is an interesting problem whether a generic nonhyperelliptic surface
in its moduli space has no automorphisms except the identity.
For surfaces in $\Bbb P^3$,
the situation is simple;
we can modify the proof of \cite{Ha2, Ex. 5.7(c)} to prove
the following.

 \proclaim{Proposition 8.4}
 If  $d\ge 5$, then smooth  surfaces in  $\Bbb C\Bbb P^3$
of degree $d$ having general moduli
 has no automorphisms except the identity.
Here {\it of general moduli} means that there is a countable
union $V$ of subvarieties of the space $\Bbb P^N$
of surfaces of degree $d$ in $\Bbb P^3$,
such that the statement $\o{Aut}(S)=1$ holds for $S\in \Bbb P^N-V$.
\qed \endproclaim

\end